\documentclass[12pt, fleqn]{article} 

\usepackage{amsmath,amssymb,amsthm}

\usepackage[numbers]{natbib} 
\usepackage[margin=0.75in,centering]{geometry}
\usepackage{graphicx,float,pdflscape}
\usepackage{enumitem}
\usepackage[utf8]{inputenc}
\usepackage[section]{placeins}

\usepackage{makecell,longtable,booktabs,array}
\usepackage[table]{xcolor} 
\definecolor{lightgray}{rgb}{0.9, 0.9, 0.9}

\usepackage{microtype} 
\usepackage[colorlinks=true, linkcolor=blue, citecolor=blue, urlcolor=blue]{hyperref}
\usepackage[capitalize,noabbrev]{cleveref}

\newtheorem{thm}{Theorem}[section]

\newtheorem{lem}{Lemma}[section]
\newtheorem{prop}{Proposition}[section]
\theoremstyle{definition} 
\newtheorem{defn}{Definition}[section]
\newtheorem{example}{Example}[section]
\theoremstyle{remark}

\usepackage{newtxtext,newtxmath}

\title{Globally Simple Heffter Arrays $H(n;k)$ with $k \equiv 1 \pmod{4}$ 
       \thanks{This research was funded by TÜBİTAK, Grant/Award Number: 124F360}}
\author{Erik Pelttari \thanks{Department of Mathematics, Koç University. Email: epelttari22@ku.edu.tr}
        \and Selda Küçükçifçi \thanks{Department of Mathematics, Koç University. Email: skucukcifci@ku.edu.tr}
        \and E. Şule Yazıcı \thanks{Department of Mathematics, Koç University. Email: eyazici@ku.edu.tr}
        }

\date{\today}

\begin{document}

\maketitle 

\begin{abstract}
 Heffter arrays are combinatorial structures used to construct orthogonal cyclic cycle decompositions and biembeddings of complete graphs onto surfaces. A Heffter array $H(m,n;s,k)$ is an $m \times n$ partially filled array with distinct nonzero entries from $\mathbb{Z}_{2nk+1}$ such that each row contains $s$ filled cells, each column contains $k$ filled cells, the elements in the filled cells form a half-set of $\mathbb{Z}_{2nk+1}$, and every row and column sums to zero modulo $2nk+1$. If these row and column sums equal zero over the integers, the structure is called an integer Heffter array. Furthermore, such an array is called globally simple if the partial sums of the entries in each row and column, evaluated in their natural order, are distinct modulo $2nk+1$. When $m=n$ and $s=k$, the array is square and denoted by $H(n;k)$. While the existence of globally simple square Heffter arrays has been established for several congruence classes, the cases where $k \equiv 1,2 \pmod{4}$ for $k > 10$ have remained an open problem \cite{dinitz2024}. In this work, we address this gap in the literature by explicitly constructing globally simple integer Heffter arrays $H(n;k)$ for the previously open cases where $k \equiv 1 \pmod{4}$ and $n \equiv 0,3 \pmod{4}$. Consequently, these constructions guarantee the existence of orthogonal cyclic $k$-cycle decompositions of the complete graph $K_{2nk+1}$ for these parameters.

\vspace{1em}
\noindent\textbf{Keywords and MSC codes:} Heffter arrays, globally simple, compatible orderings, cycle decompositions. 05B30, 05C51
\end{abstract}

\section{Introduction}

\subsection{Overview and Definitions}

Heffter arrays were first introduced by Archdeacon \cite{archdeacon2015} in 2015 as a combinatorial tool to construct orthogonal cyclic cycle decompositions and biembeddings of complete graphs onto surfaces. As their name suggests, these structures can be viewed as a two-dimensional generalization of Heffter's first difference problem, originally posed in 1896 \cite{heffter1896}. While the classical difference problem seeks to partition the set $\{1,2,\dots,3m\}$ into triples that sum to zero modulo $6m+1$, Heffter arrays extend this concept to a matrix format where the entries in every row and column must sum to zero modulo $2nk+1$.

\begin{defn}\label{def:heffter}
A \textbf{Heffter array} $H(m,n; s,k)$ is an $m \times n$ matrix with entries from $\mathbb{Z}_{2nk+1}$, where $ms = nk$, such that:
\begin{enumerate}[label=(\roman*)]
    \item Each row contains exactly $s$ filled cells and each column contains exactly $k$ filled cells.
    \item For every $x \in \mathbb{Z}_{2nk+1} \setminus \{0\}$, exactly one of $x$ or $-x$ appears in the array.
    \item The sum of the elements in every row and every column is $0 \pmod{2nk+1}$.
\end{enumerate}
When the array is square (i.e., $m=n$ and thus $s=k$), it is denoted simply as $H(n;k)$. 
\end{defn}

\begin{example}
\label{ex:H86}
An $H(8; 6)$ where the entries come from $\mathbb{Z}_{97}$.

\begin{center}
\renewcommand{\arraystretch}{1.25} 
\begin{tabular}{|c|c|c|c|c|c|c|c|}
\hline
$-1$ & $5$ & $2$ & $-7$ & $-9$ & $10$ & & \\ \hline
$3$ & $-4$ & $-6$ & $8$ & $11$ & $-12$ & & \\ \hline
& & $-13$ & $17$ & $14$ & $-19$ & $-21$ & $22$ \\ \hline
& & $15$ & $-16$ & $-18$ & $20$ & $23$ & $-24$ \\ \hline
$-33$ & $34$ & & & $-25$ & $29$ & $26$ & $-31$ \\ \hline
$35$ & $-36$ & & & $27$ & $-28$ & $-30$ & $32$ \\ \hline
$38$ & $-43$ & $-45$ & $46$ & & & $-37$ & $41$ \\ \hline
$-42$ & $44$ & $47$ & $-48$ & & & $39$ & $-40$ \\ \hline
\end{tabular}
\end{center}
\end{example}

\begin{defn}\label{def:tight}
A Heffter array $H(m,n;n,m)$ with no empty cells is called a \textbf{tight Heffter array} and is denoted $H(m,n)$.
\end{defn}

\begin{defn}\label{def:integer}
If the elements in every row and column sum to zero over the integers $\mathbb{Z}$ (rather than merely modulo $2nk+1$), the array is called an \textbf{integer Heffter array}.
\end{defn}

\begin{defn}
The \textbf{support} of a partially filled array is the set of the absolute values of its non-empty entries. By condition (ii) of \cref{def:heffter}, it follows that the support of a Heffter array is precisely the set $\{1, 2, \dots, nk\}$. 
\end{defn}

While Heffter arrays are interesting combinatorial objects in their own right, their utility lies in their applications. Specifically, these arrays, with several additional properties, provide a recipe for decomposing complete graphs \cite{archdeacon2015, costa2018} and embedding them onto surfaces \cite{archdeacon2015, cavenagh2020}. The next set of definitions covers these graph-theoretic and topological concepts.

\begin{defn}\label{def:orthogonal}
A cycle decomposition of a graph $G$ is a set of cycles, all of the same length, whose union contains every edge of $G$ exactly once. Two such decompositions, $\mathcal{C}_1$ and $\mathcal{C}_2$, are \textbf{orthogonal} if any cycle from $\mathcal{C}_1$ and any cycle from $\mathcal{C}_2$ intersect in at most one edge \cite{caro1999}.
\end{defn}

\begin{defn}\label{def:simple}
Let $T$ be a finite subset of an additive group $G$. Given an ordering $\omega = (t_1, t_2, \dots, t_k)$ of the elements in $T$, let $\displaystyle{s_i = \sum_{j=1}^i t_j}$ be the $i$-th partial sum of $\omega$. The ordering $\omega$ is said to be \textbf{simple} if its partial sums $s_1, s_2, \dots, s_k$ are pairwise distinct. Note that $\omega = (t_1, \ \dots, \ t_k)$ is simple if and only if any cyclic shift $\omega^i = (t_i, \ \dots, \ t_k, \ t_1, \ \dots, \ t_{i-1})$ is simple.   
\end{defn}

\begin{defn}\label{def:globsimp}
In the context of a Heffter array $H(m,n;s,k)$, the filled entries of any row or column form a finite subset of the additive group $\mathbb{Z}_{2nk+1}$. A Heffter array is said to be \textbf{simple} if there exists a simple ordering for the elements of every row and every column. Given an arbitrary Heffter array, determining whether such simple orderings exist can be computationally difficult. Because of this, arrays that are constructed to be simple are of particular interest. A Heffter array is defined to be \textbf{globally simple} if all of its rows and columns are simple under their natural left-to-right and top-to-bottom orderings, respectively \cite{costa2018}. If a Heffter array is simple, the rows and columns can be used as difference sets to construct cyclic cycle decompositions of $K_{2nk+1}$. Specifically, the rows can be used to produce a cyclic decomposition of the complete graph $K_{2nk+1}$ into $s$-cycles, the columns can be used to produce a cyclic decomposition of $K_{2nk+1}$ into $k$-cycles, and these two decompositions are orthogonal to each other. 
\end{defn}

\begin{example}\label{ex:nonsimprow}
Consider the row entries $R_1 = \{-1, 5, 2, -7, -9, 10\}$ from the Heffter array $H(8;6)$ in \cref{ex:H86} with entries in $\mathbb{Z}_{97}$. We evaluate two different orderings of these elements:

\begin{enumerate}\item \textbf{The Natural Ordering:} Let $\omega_1 = (-1, 5, 2, -7, -9, 10)$. The sequence of partial sums is $S = (-1, 4, 6, -1, -10, 0)$. In the context of $K_{97}$, this ordering corresponds to a walk starting at vertex $0$. Because the partial sum $-1$ occurs twice (at $s_1$ and $s_4$), the walk visits vertex $-1$ twice before returning to $0$. This produces two disjoint $3$-cycles: one on vertices $\{-1, 4, 6\}$ and another on $\{0, -1, -10\}$. Thus, $\omega_1$ is not a simple ordering.\item \textbf{A Simple Reordering:} Let $\omega_1^* = (-1, 5, 2, -9, -7, 10)$, obtained by transposing the fourth and fifth elements. The resulting partial sums are $S^* = (-1, 4, 6, -3, -10, 0)$. Since the elements of $S^*$ are pairwise distinct and non-zero, $\omega_1^*$ is a simple ordering. This ordering successfully traces a single $6$-cycle in $K_{97}$ with vertices $(0, -1, 4, 6, -3, -10)$.
\end{enumerate}

In this specific case, finding a simple reordering was straightforward, but this is not true in general; determining whether an arbitrary set of integers admits a simple ordering is a non-trivial problem \cite{archdeacon2016}. Because one of our goals is to use Heffter arrays to produce orthogonal cycle decompositions, we prefer constructions that are globally simple, as they provide these cycle decompositions directly without the need for additional reordering.
\end{example}

\begin{defn}\label{def:compatible}
The cyclic ordering of the entries in a row or column of a partially filled array (such as a Heffter array) can be viewed as a permutation on the set of entries in filled cells of the array. Let $\omega_r$ denote the permutation formed by the composition of the row cycles, and let $\omega_c$ denote the corresponding permutation formed by composition of the column cycles. If the composition $\omega_r \circ \omega_c$ is a single cycle of length equal to the total number of filled cells in the array, we say the row and column orderings, $\omega_r$ and $\omega_c$ respectively, are \textbf{compatible}.
\end{defn}

\begin{defn}\label{def:orientable}
A surface is \textbf{orientable} if it admits a consistent local orientation around every point. Informally, it is a surface that does not contain a Möbius band. An \textbf{orientable embedding} of a graph is an embedding of the graph onto an orientable surface such that edges intersect only at their incident vertices, and every face is a 2-cell (homeomorphic to an open disk).
\end{defn}

\begin{defn}\label{def:biembedding}
A cellular embedding of a graph $G$ on a surface is called a \textbf{biembedding} if its faces are $2$-colorable. In a biembedding, every edge of $G$ is on exactly one face of each color. If the faces of the biembedding are simple cycles, the biembedding is said to be simple. 
\end{defn}

In our context, given a Heffter array $H(m,n;s,k)$, we would like to determine if we can produce orthogonal cycle decompositions $\mathcal{C}_1$ and $\mathcal{C}_2$ of the complete graph $K_{2nk+1}$ such that the cycles in $\mathcal{C}_1$ have length $s$ and the cycles in $\mathcal{C}_2$ have length $k$. In addition, we would like to determine when we can construct a biembedding of $K_{2nk+1}$  such that the cycles in $\mathcal{C}_1$ correspond to the faces in one color class and the cycles in $\mathcal{C}_2$ correspond to the faces in the second color class. 

This connection between the structure of Heffter arrays and biembeddings was established by Archdeacon:

\begin{thm}\textup{\cite{archdeacon2015}}
\label{thm:archdeacon}
Given a Heffter array $H(m,n; s,k)$ with compatible orderings $\omega_r$ on the rows and $\omega_c$ on the columns, there exists an orientable embedding of $K_{2nk+1}$ such that every edge is on a face of size $s$ and a face of size $k$. Moreover, if $\omega_r$ and $\omega_c$ are both simple, then all faces are simple cycles.
\end{thm}

If a Heffter array possesses compatible orderings on the rows and columns but fails to be simple, we still obtain a biembedding on an orientable surface; however, the rows or columns (whichever does not have a simple ordering) will produce circuits instead of cycles, resulting in faces that are not simple. Conversely, if the array is simple but lacks a compatible ordering, its rows and columns still produce orthogonal cyclic decompositions of the complete graph. Topologically, however, the lack of a single full-length permutation cycle causes the structure to pinch at one or more of the vertices, resulting in an embedding on a pseudosurface rather than on an orientable surface \cite{archdeacon2015}.

The study of Heffter arrays naturally begins with the question of their existence. Square Heffter arrays trivially require $3 \le k \le n$, and for square integer Heffter arrays, we additionally require $nk \equiv 0, 3 \pmod{4}$ \cite{archdeacon2_2015}. Through the combined efforts of several authors, the existence problem for these arrays was completely solved:

\begin{thm}
\textup{\cite{archdeacon2_2015, cavenagh2019, dinitz2017}} 
A square Heffter array $H(n;k)$ exists if and only if $3 \le k \le n$, and a square integer Heffter array $H(n;k)$ exists if and only if $3 \le k \le n$ and $nk \equiv 0,3 \pmod{4}$.
\end{thm}

With existence established, the focus shifted toward finding Heffter arrays that also possess the \textbf{globally simple} property, as these arrays are guaranteed to yield orthogonal cyclic decompositions of complete graphs \cite{kucukcifci}. The baseline existence of globally simple arrays for smaller values of $k$ was established by Costa, Morini, Pasotti, and Pellegrini:

\begin{thm}\textup{\cite{costa2018}}
\label{thm:costa_small_k}
For $3 \le k \le 10$, a globally simple Heffter array $H(n;k)$ exists if and only if $k \le n$ and $nk \equiv 0,3 \pmod{4}$.
\end{thm}

For broader classes where $k \equiv 0,3 \pmod{4}$, Burrage et al.\ constructed globally simple Heffter arrays under specific conditions:

\begin{thm}
\label{thm:burrage_classes}
\textup{\cite{burrage2020}}
There exist globally simple Heffter arrays $H(n;k)$ for $3 \le k \le n$ in the following cases:
\begin{enumerate}
    \item $k \equiv 0 \pmod{4}$,
    \item $n \equiv 1 \pmod{4}$ and $k \equiv 3 \pmod{4}$, and
    \item $n \equiv 0 \pmod{4}$ and $k \equiv 3 \pmod{4}$, provided $n \gg k$.
\end{enumerate}
\end{thm}

However, the existence of globally simple Heffter arrays for the cases where $k \equiv 1,2 \pmod{4}$ and $k > 10$ has remained an open problem. In this paper, we address the problem of globally simple Heffter arrays for $k \equiv 1 \pmod{4}$ with $k \ge 13$ by stating and proving the following results:

\begin{thm}
\label{thm:n3k1}
For $n \equiv 3 \pmod{4}$, $k \equiv 1 \pmod{4}$ and $n > k \ge 13$, there exists a globally simple integer Heffter array $H(n;k)$.
\end{thm}

\begin{thm}
\label{thm:n0k1}
For $n \equiv 0 \pmod{4}$, $k \equiv 1 \pmod{4}$ and $n > k \ge 13$, there exists a globally simple integer Heffter array $H(n;k)$.
\end{thm}

Regarding compatible orderings, which imply biembeddings on orientable surfaces, the following theorem states the necessary conditions on a general Heffter array $H(m,n;s,k)$.

\begin{thm}
\label{thm:compatible_ness}
\textup{\cite{cavenagh2020, costa2020, dinitz2_2017, dinitz2024}}
If there exist compatible orderings $\omega_r$ and $\omega_c$ for a Heffter array $H(m,n;s,k)$, then either: $m,n,s$ and $k$ are all odd; $m$ is odd, $n$ is even and $s$ is even; or $m$ is even, $n$ is odd and $k$ is even.
\end{thm}

For square integer Heffter arrays, this implies that $nk$ is odd. This requirement leaves two primary cases for odd $n$ and $k$. In their work on biembeddings of cycle systems, Cavenagh et al. proved the existence of compatible orderings for $n \equiv 1 \pmod{4}$ and $k \equiv 3 \pmod{4}$ with sporadic exceptions \cite{cavenagh2020}. 

In the case where $n \equiv 3 \pmod{4}$ and $k \equiv 1 \pmod{4}$, a proposition from Costa et~al. shows that when $\gcd(n,k-2) = 1$, Heffter arrays constructed in a specific manner have compatible orderings.

\begin{prop}
\textup{\cite{costa2020}}
  Let $3 \le k < n$ be odd integers. Let the entries of $H(n;k)$ be entirely contained in $k$ diagonals with two sets of empty diagonals, both of which have width $\frac{n-k}{2}$. If $\gcd(n,k-2) = 1$, then there exist compatible orderings of the rows and columns of $H(n;k)$.  
\end{prop}

The Heffter array constructed in \Cref{sec:construction31} when $n \equiv 3 \pmod{4}$ and $k \equiv 1 \pmod{4}$ have entries entirely contained in $k$ diagonals, with exactly two sets of empty diagonals of width $\frac{n-k}{2}$. As a result, this proposition immediately yields the following theorem.

\begin{thm}
Let $n \equiv 3 \pmod{4}$, $k \equiv 1 \pmod{4}$, $n > k \ge 13$, and $\gcd(n,k-2) = 1$. Then there exists a biembedding of $K_{2nk+1}$ on an orientable surface such that each face is a simple cycle of length $k$.  
\end{thm}

The remainder of this paper is organized as follows. In \cref{sec:construction31}, we describe the construction of the array $H(n;k)$ for $n \equiv 3 \pmod{4}$ and $k \equiv 1 \pmod{4}$. In \cref{sec:proof_simple31}, we prove that this construction yields a globally simple array by showing that the row and column partial sums are distinct modulo $2nk+1$. In \cref{sec:construction01} and \cref{sec:proof_simple01}, we show that a globally simple $H(n;k)$ exists when $n \equiv 0 \pmod{4}$.


\section{$H(n;k)$ Construction for $n \equiv 3 \pmod{4}$ and $k \equiv 1 \pmod{4}$} \label{sec:construction31}

We construct integer globally simple Heffter arrays $H(n;k)$ for all $n \equiv 3 \pmod{4}$ and $k \equiv 1 \pmod{4}$ with $n > k \geq 13$. The construction proceeds by initializing a base $H(n;5)$ array, as described by Dinitz and Wanless \cite{dinitz2017}, and augmenting with sets of 4 diagonals until $k$ is reached.

\subsection{Diagonal Definition}

We define a diagonal of length $\ell$, denoted by $\text{diag}(r, c, v, \delta, \delta_v, \ell)$, as the set of $\ell$ filled cells starting at cell $(r, c)$ with value $v$. Subsequent cells are determined by incrementing the row and column indices by $\delta$ and the value by $\delta_v$. Formally, for $j \in \{0, 1, \dots, \ell-1\}$, the diagonal contains the entry $v + j\delta_v$ in cell $(r + j\delta, c + j\delta)$. Unless otherwise stated, all rows and columns are indexed by $\{1,2,\dots,n\}$ and operations are performed modulo $n$.

\subsection{The Base Array $H(n;5)$ when $n \equiv 3 \pmod{4}$}

Let $h=(n+1)/2$ and $q=(n-3)/4$. The base $H(n;5)$ is constructed using the following diagonals and ad hoc entries.

\paragraph{Diagonal Assignments}
\begin{align*}
    &\text{A: diag}(h+1, h+1, (h-2), 1, -1, h-3) & &\text{H: diag}(h+2, h+1, -(5n-q-3), 2, 1, q) \\
    &\text{B: diag}(3, 3, -(n-3), 1, 1, h-3) & &\text{I: diag}(h, h+1, (4n+q+1), 2, 1, q)\\
    &\text{C: diag}(2, 3, (4n), 2, -1, q) & &\text{J: diag}(h+1, h+2, (3n+q+4), 2, 1, q)\\
    &\text{D: diag}(3, 2, -(4n+1), 2, -1, q) & &\text{K: diag}(h+1, 2, -(n+2), 1, -2, h-2)\\
    &\text{E: diag}(3, 4, (5n-3), 2, -1, q) & &\text{L: diag}(h+1, 1, (n+1), 1, 2, h-1)\\
    &\text{F: diag}(4, 3, -(3n+4), 2, -1, q) & &\text{M: diag}(1, h+1, -(3n), 1, 2, h-1)\\
    &\text{G: diag}(h+1, h, -(4n-q), 2, 1, q) & &\text{N: diag}(2, h+1, (3n-1), 1, -2, h-2)
\end{align*}

\paragraph{Ad Hoc Cell Values}
The following specific entries are set directly:

\begin{align*}
    &H[1,1] = n-1, &&H[1,2] = 5n-2, &&H[1,h] = 2n+2, &&H[1,n] = -5n+1, \\
    &H[2,1] = -3n-3, &&H[2,2] = -n+2, \\
    &H[h,1] = -2n+1, &&H[h,h] = -n, &&H[h,n] = 2n+1, \\
    &H[n-1,n-1] = -h+1, &&H[n-1,n] = 5n, \\
    &H[n,1] = 3n+2, &&H[n,h] = -2n, &&H[n,n-1] = -3n-1, &&H[n,n] = 1. \\
\end{align*}

\subsection{Augmentation with $k-5$ New Diagonals}

To increase the number of filled cells in each row and column from 5 to $k$, we augment the base construction with $t = (k-5)/4$ sets of four diagonals. For each $i \in \{1, 2, \dots, t\}$, we define the following four diagonals:

\begin{equation*}
\begin{array}{lll}
\toprule
\textbf{Diagonal} & \textbf{Definition} & \textbf{Support} \\
\midrule
\mathbf{D_{1+2i}} & \text{diag}(1+2i, 1, (n + 4ni + 1), 1, 2, n) & \{n+4ni+1, n+4ni+3, \dots, 3n+4ni-1\} \\
\addlinespace
\mathbf{D_{2+2i}} & \text{diag}(2+2i, 1, (-3n - 4ni - 1), 1, -2, n) & \{3n+4ni+1, 3n+4ni+3, \dots, 5n+4ni-1\} \\
\addlinespace
\mathbf{D_{h+2i}} & \text{diag}(1+2i, h+1, (-n - 4ni - 2), 1, -2, n) & \{n+4ni+2, n+4ni+4, \dots, 3n+4ni\} \\
\addlinespace
\mathbf{D_{h+1+2i}} & \text{diag}(2+2i, h+1, (3n + 4ni + 2), 1, 2, n) & \{3n+4ni+2, 3n+4ni+4, \dots, 5n+4ni\} \\
\bottomrule
\end{array}
\end{equation*}

Combined with the base support $\{1, \dots, 5n\}$ from the $H(n;5)$, the constructed $H(n;k)$ has the required support $\{1, 2, \dots, nk\}$.

\subsection{Row and Column Sums}

The rows and columns of the base $H(n;5)$ sum to zero. We verify the above diagonals sum to zero as well. Consider the contribution of the $i$-th set of diagonals to any fixed row $r$. The row $r$ intersects $\mathbf{D_{1+2i}}$ at some column $c$ with value $v$, and intersects $\mathbf{D_{h+2i}}$ at column $c'$ with value $v'$. By the diagonal definitions:
\[ v + v' = (n+4ni+1+2((r-1-2i) \bmod{n})) + (-n-4ni-2-2((r-1-2i) \bmod{n})) = -1. \]
Similarly, the values from $\mathbf{D_{2+2i}}$ and $\mathbf{D_{h+1+2i}}$ sum to $+1$. Thus, the total contribution of these diagonals to row $r$ is $-1 + 1 = 0$. For any fixed column $c$, the values from $\mathbf{D_{1+2i}}$ and $\mathbf{D_{2+2i}}$ sum to $-2n$, while the values from $\mathbf{D_{h+2i}}$ and $\mathbf{D_{h+1+2i}}$ sum to $+2n$. Thus, the total contribution of these diagonals to column $c$ is $-2n + 2n = 0$.

Since the support is correct and every row and column sum is zero over the integers, this construction produces a valid integer Heffter array $H(n;k)$.

\subsection{Example: Globally Simple $H(19;17)$}

The 5 light gray diagonals in \Cref{tab:H1917} are the base $H(n;5)$. The starting cells for the 12 additional diagonals are highlighted in dark gray.

\begin{table}[h]
    \centering
    \caption{Globally Simple Heffter Array $H(19;17)$}
    \label{tab:H1917}
    \setlength{\tabcolsep}{4pt}
    \scriptsize
        \begin{tabular}{|*{19}{c|}} 
            \hline
\cellcolor{black!15}18 & \cellcolor{black!15}93 &  & 311 & -275 & 239 & -203 & 167 & -131 & \cellcolor{black!15}40 & \cellcolor{black!15}-57 &  & -310 & 274 & -238 & 202 & -166 & 130 & \cellcolor{black!15}-94 \\
\hline\cellcolor{black!15}-60 & \cellcolor{black!15}-17 & \cellcolor{black!15}76 &  & 313 & -277 & 241 & -205 & 169 & -133 & \cellcolor{black!15}56 & \cellcolor{black!15}-55 &  & -312 & 276 & -240 & 204 & -168 & 132 \\
\hline \cellcolor{black!40}96 & \cellcolor{black!15}-77 & \cellcolor{black!15}-16 & \cellcolor{black!15}92 &  & 315 & -279 & 243 & -207 & 171 & \cellcolor{black!40}-97 & \cellcolor{black!15}54 & \cellcolor{black!15}-53 &  & -314 & 278 & -242 & 206 & -170 \\
\hline \cellcolor{black!40}-134 & 98 & \cellcolor{black!15}-61 & \cellcolor{black!15}-15 & \cellcolor{black!15}75 &  & 317 & -281 & 245 & -209 & \cellcolor{black!40}135 & -99 & \cellcolor{black!15}52 & \cellcolor{black!15}-51 &  & -316 & 280 & -244 & 208 \\
\hline \cellcolor{black!40}172 & -136 & 100 & \cellcolor{black!15}-78 & \cellcolor{black!15}-14 & \cellcolor{black!15}91 &  & 319 & -283 & 247 & \cellcolor{black!40}-173 & 137 & -101 & \cellcolor{black!15}50 & \cellcolor{black!15}-49 &  & -318 & 282 & -246 \\
\hline \cellcolor{black!40}-210 & 174 & -138 & 102 & \cellcolor{black!15}-62 & \cellcolor{black!15}-13 & \cellcolor{black!15}74 &  & 321 & -285 & \cellcolor{black!40}211 & -175 & 139 & -103 & \cellcolor{black!15}48 & \cellcolor{black!15}-47 &  & -320 & 284 \\
\hline \cellcolor{black!40}248 & -212 & 176 & -140 & 104 & \cellcolor{black!15}-79 & \cellcolor{black!15}-12 & \cellcolor{black!15}90 &  & 323 & \cellcolor{black!40}-249 & 213 & -177 & 141 & -105 & \cellcolor{black!15}46 & \cellcolor{black!15}-45 &  & -322 \\
\hline \cellcolor{black!40}-286 & 250 & -214 & 178 & -142 & 106 & \cellcolor{black!15}-63 & \cellcolor{black!15}-11 & \cellcolor{black!15}73 &  & \cellcolor{black!40}287 & -251 & 215 & -179 & 143 & -107 & \cellcolor{black!15}44 & \cellcolor{black!15}-43 &  \\
\hline & -288 & 252 & -216 & 180 & -144 & 108 & \cellcolor{black!15}-80 & \cellcolor{black!15}-10 & \cellcolor{black!15}89 &  & 289 & -253 & 217 & -181 & 145 & -109 & \cellcolor{black!15}42 & \cellcolor{black!15}-41 \\
\hline\cellcolor{black!15}-37 &  & -290 & 254 & -218 & 182 & -146 & 110 & \cellcolor{black!15}-64 & \cellcolor{black!15}-19 & \cellcolor{black!15}81 &  & 291 & -255 & 219 & -183 & 147 & -111 & \cellcolor{black!15}39 \\
\hline\cellcolor{black!15}20 & \cellcolor{black!15}-21 &  & -292 & 256 & -220 & 184 & -148 & 112 & \cellcolor{black!15}-72 & \cellcolor{black!15}8 & \cellcolor{black!15}65 &  & 293 & -257 & 221 & -185 & 149 & -113 \\
\hline-115 & \cellcolor{black!15}22 & \cellcolor{black!15}-23 &  & -294 & 258 & -222 & 186 & -150 & 114 & \cellcolor{black!15}-88 & \cellcolor{black!15}7 & \cellcolor{black!15}82 &  & 295 & -259 & 223 & -187 & 151 \\
\hline153 & -117 & \cellcolor{black!15}24 & \cellcolor{black!15}-25 &  & -296 & 260 & -224 & 188 & -152 & 116 & \cellcolor{black!15}-71 & \cellcolor{black!15}6 & \cellcolor{black!15}66 &  & 297 & -261 & 225 & -189 \\
\hline-191 & 155 & -119 & \cellcolor{black!15}26 & \cellcolor{black!15}-27 &  & -298 & 262 & -226 & 190 & -154 & 118 & \cellcolor{black!15}-87 & \cellcolor{black!15}5 & \cellcolor{black!15}83 &  & 299 & -263 & 227 \\
\hline229 & -193 & 157 & -121 & \cellcolor{black!15}28 & \cellcolor{black!15}-29 &  & -300 & 264 & -228 & 192 & -156 & 120 & \cellcolor{black!15}-70 & \cellcolor{black!15}4 & \cellcolor{black!15}67 &  & 301 & -265 \\
\hline-267 & 231 & -195 & 159 & -123 & \cellcolor{black!15}30 & \cellcolor{black!15}-31 &  & -302 & 266 & -230 & 194 & -158 & 122 & \cellcolor{black!15}-86 & \cellcolor{black!15}3 & \cellcolor{black!15}84 &  & 303 \\
\hline305 & -269 & 233 & -197 & 161 & -125 & \cellcolor{black!15}32 & \cellcolor{black!15}-33 &  & -304 & 268 & -232 & 196 & -160 & 124 & \cellcolor{black!15}-69 & \cellcolor{black!15}2 & \cellcolor{black!15}68 &  \\
\hline & 307 & -271 & 235 & -199 & 163 & -127 & \cellcolor{black!15}34 & \cellcolor{black!15}-35 &  & -306 & 270 & -234 & 198 & -162 & 126 & \cellcolor{black!15}-85 & \cellcolor{black!15}-9 & \cellcolor{black!15}95 \\
\hline\cellcolor{black!15}59 &  & 309 & -273 & 237 & -201 & 165 & -129 & \cellcolor{black!15}36 & \cellcolor{black!15}-38 &  & -308 & 272 & -236 & 200 & -164 & 128 & \cellcolor{black!15}-58 & \cellcolor{black!15}1 \\
\hline
\end{tabular}
\end{table}

\section{$H(n;k)$ is globally simple for $n \equiv 3 \pmod{4}$ and $k \equiv 1 \pmod{4}$} \label{sec:proof_simple31}

In this section, we prove \Cref{thm:n3k1} by showing the integer Heffter arrays constructed in \Cref{sec:construction31} are globally simple. 

\begin{defn}\label{def:partialsums}
Let $S_C(c, p)$ and $S_R(r, p)$ denote the $p$-th partial sums of column $c$ and row $r$, respectively, where $1 \le c, r \le n$ and $1 \le p \le k$.
\end{defn}

\subsection{Columns are simple}

To prove a column is simple, we must compute its partial sums and ensure they are pairwise distinct. For each column $c \in \{1,2,\dots,n\}$, the column partial sums $S_C(c, p)$ are calculated top-to-bottom, starting in cell $(c+2, c)$. As mentioned in \Cref{def:simple}, an ordering is simple if and only if any cyclic shift of that ordering is simple. Therefore, our choice of starting cell does not affect the simplicity proof. 

\subsubsection{Example: Globally Simple $H(19;17)$ Column Partial Sums}
We illustrate the partial sum computation using a globally simple $H(19;17)$. \Cref{tab:H1917colsums} displays the computed partial sums for each column. For instance, the values in the first row of the table correspond to the sequence $S_C(1, 1), S_C(1, 2), \dots, S_C(1, 17)$. Shaded cells indicate partial sums derived from the base array $H(19;5)$.

\begin{table}[h]
    \centering
    \caption{Column Partial Sums Modulo $2nk+1=647$ for the Globally Simple $H(19;17)$ }
    \label{tab:H1917colsums}
    \setlength{\tabcolsep}{4pt} 
    \scriptsize
        
    \begin{tabular}{lc*{17}{c}} 
    \toprule
    Col 1 && 96 & $-$38 $\equiv$ 609 & 134 & 571 & 172 & 533 & \cellcolor{lightgray}496 & \cellcolor{lightgray}516 & 401 & 554 & 363 & 592 & 325 & 630 & \cellcolor{lightgray}42 & \cellcolor{lightgray}60 & \cellcolor{lightgray}0 \\ 
    Col 2 && 98 & 609 & 136 & 571 & 174 & 533 & \cellcolor{lightgray}512 & \cellcolor{lightgray}534 & 417 & 572 & 379 & 610 & 341 & 1 & \cellcolor{lightgray}94 & \cellcolor{lightgray}77 & \cellcolor{lightgray}0 \\ 
    Col 3 && 100 & 609 & 138 & 571 & 176 & 533 & \cellcolor{lightgray}510 & \cellcolor{lightgray}534 & 415 & 572 & 377 & 610 & 339 & 1 & \cellcolor{lightgray}77 & \cellcolor{lightgray}61 & \cellcolor{lightgray}0 \\ 
    Col 4 && 102 & 609 & 140 & 571 & 178 & 533 & \cellcolor{lightgray}508 & \cellcolor{lightgray}534 & 413 & 572 & 375 & 610 & 337 & 1 & \cellcolor{lightgray}93 & \cellcolor{lightgray}78 & \cellcolor{lightgray}0 \\ 
    Col 5 && 104 & 609 & 142 & 571 & 180 & 533 & \cellcolor{lightgray}506 & \cellcolor{lightgray}534 & 411 & 572 & 373 & 610 & 335 & 1 & \cellcolor{lightgray}76 & \cellcolor{lightgray}62 & \cellcolor{lightgray}0 \\ 
    Col 6 && 106 & 609 & 144 & 571 & 182 & 533 & \cellcolor{lightgray}504 & \cellcolor{lightgray}534 & 409 & 572 & 371 & 610 & 333 & 1 & \cellcolor{lightgray}92 & \cellcolor{lightgray}79 & \cellcolor{lightgray}0 \\ 
    Col 7 && 108 & 609 & 146 & 571 & 184 & 533 & \cellcolor{lightgray}502 & \cellcolor{lightgray}534 & 407 & 572 & 369 & 610 & 331 & 1 & \cellcolor{lightgray}75 & \cellcolor{lightgray}63 & \cellcolor{lightgray}0 \\ 
    Col 8 && 110 & 609 & 148 & 571 & 186 & 533 & \cellcolor{lightgray}500 & \cellcolor{lightgray}534 & 405 & 572 & 367 & 610 & 329 & 1 & \cellcolor{lightgray}91 & \cellcolor{lightgray}80 & \cellcolor{lightgray}0 \\ 
    Col 9 && 112 & 609 & 150 & 571 & 188 & 533 & \cellcolor{lightgray}498 & \cellcolor{lightgray}534 & 403 & 572 & 365 & 610 & 327 & 1 & \cellcolor{lightgray}74 & \cellcolor{lightgray}64 & \cellcolor{lightgray}0 \\ 
    Col 10 && 114 & 609 & 152 & 571 & 190 & 533 & \cellcolor{lightgray}495 & \cellcolor{lightgray}535 & 402 & 573 & 364 & 611 & 326 & 2 & \cellcolor{lightgray}91 & \cellcolor{lightgray}72 & \cellcolor{lightgray}0 \\ 
    Col 11 && 116 & 609 & 154 & 571 & 192 & 533 & \cellcolor{lightgray}476 & \cellcolor{lightgray}532 & 435 & 570 & 397 & 608 & 359 & 646 & \cellcolor{lightgray}80 & \cellcolor{lightgray}88 & \cellcolor{lightgray}0 \\ 
    Col 12 && 118 & 609 & 156 & 571 & 194 & 533 & \cellcolor{lightgray}478 & \cellcolor{lightgray}532 & 433 & 570 & 395 & 608 & 357 & 646 & \cellcolor{lightgray}64 & \cellcolor{lightgray}71 & \cellcolor{lightgray}0 \\ 
    Col 13 && 120 & 609 & 158 & 571 & 196 & 533 & \cellcolor{lightgray}480 & \cellcolor{lightgray}532 & 431 & 570 & 393 & 608 & 355 & 646 & \cellcolor{lightgray}81 & \cellcolor{lightgray}87 & \cellcolor{lightgray}0 \\ 
    Col 14 && 122 & 609 & 160 & 571 & 198 & 533 & \cellcolor{lightgray}482 & \cellcolor{lightgray}532 & 429 & 570 & 391 & 608 & 353 & 646 & \cellcolor{lightgray}65 & \cellcolor{lightgray}70 & \cellcolor{lightgray}0 \\ 
    Col 15 && 124 & 609 & 162 & 571 & 200 & 533 & \cellcolor{lightgray}484 & \cellcolor{lightgray}532 & 427 & 570 & 389 & 608 & 351 & 646 & \cellcolor{lightgray}82 & \cellcolor{lightgray}86 & \cellcolor{lightgray}0 \\ 
    Col 16 && 126 & 609 & 164 & 571 & 202 & 533 & \cellcolor{lightgray}486 & \cellcolor{lightgray}532 & 425 & 570 & 387 & 608 & 349 & 646 & \cellcolor{lightgray}66 & \cellcolor{lightgray}69 & \cellcolor{lightgray}0 \\ 
    Col 17 && 128 & 609 & 166 & 571 & 204 & 533 & \cellcolor{lightgray}488 & \cellcolor{lightgray}532 & 423 & 570 & 385 & 608 & 347 & 646 & \cellcolor{lightgray}83 & \cellcolor{lightgray}85 & \cellcolor{lightgray}0 \\ 
    Col 18 && 130 & 609 & 168 & 571 & 206 & 533 & \cellcolor{lightgray}490 & \cellcolor{lightgray}532 & 421 & 570 & 383 & 608 & 345 & 646 & \cellcolor{lightgray}67 & \cellcolor{lightgray}58 & \cellcolor{lightgray}0 \\ 
    Col 19 && 132 & 609 & 170 & 571 & 208 & 533 & \cellcolor{lightgray}492 & \cellcolor{lightgray}531 & 418 & 569 & 380 & 607 & 342 & 645 & \cellcolor{lightgray}93 & \cellcolor{lightgray}94 & \cellcolor{lightgray}0 \\ 
    \bottomrule
    \end{tabular}

\end{table}

\subsubsection{Column Partial Sums Are Distinct Modulo $2nk+1$}

To prove that this construction produces a globally simple array in the general case, we first state several Lemmas used to derive the column partial sums for specific sets of columns.

\begin{lem}
Let $H(n;5)$ be the base array generated via the Dinitz and Wanless construction \cite{dinitz2017}. The non-empty entries of this array are exactly those detailed in \Cref{tab:Hn5colentries}, positioned top-to-bottom such that in each column $c$, two consecutive entries start at $(c+h-1, c)$ and three consecutive entries start at $(c-1, c)$.
\end{lem}

\begin{table}[h]
\centering
\caption{Base $H(n;5)$ Column Entries}
\label{tab:Hn5colentries}
\setlength{\tabcolsep}{4pt} 
\footnotesize
\begin{tabular}{lcccccc}
\toprule
\textbf{Column} & \multicolumn{2}{c}{\textbf{Two Consecutive Entries}} && \multicolumn{3}{c}{\textbf{Three Consecutive Entries}} \\
\midrule
Column 1 & $-2n+1$, & $n+1$, && $3n+2$, & $n-1$, & $-3n-3$ \\
\addlinespace
Column 2 & $-n-2$, & $n+3$, && $5n-2$, & $-n+2$, & $-4n-1$ \\
\midrule
Odd $3 \le c \le h-1$ & $-n-2c+2$, & $n+2c-1$, && $4n-(c-1)/2+1$, & $-n+c$, & $-3n-(c-1)/2-3$ \\
\addlinespace
Even $4 \le c \le h-2$ & $-n-2c+2$, & $n+2c-1$, && $5n-c/2-1$, & $-n+c$, & $-4n-c/2$ \\
\midrule
Column $h$ & $-2n$, & $2n+2$, && $(19n-5)/4$, & $-n$, & $-(15n+3)/4$ \\
\midrule
Odd $h+1 \le c \le n-2$ & $-4n+2c-3$, & $4n-2c+2$, && $4n+(c-1)/2$, & $n-c$, & $-5n+(c-1)/2+2$ \\
\addlinespace
Even $h+2 \le c \le n-3$ & $-4n+2c-3$, & $4n-2c+2$, && $3n+c/2+2$, & $n-c$, & $-4n+c/2-1$ \\
\midrule
Column $n-1$ & $-2n-5$, & $2n+4$, && $(7n+3)/2$, & $-(n-1)/2$, & $-3n-1$ \\
\addlinespace
Column $n$ & $-2n-3$, & $2n+1$, && $5n$, & $1$, & $-5n+1$ \\
\bottomrule
\end{tabular}
\end{table}

\begin{lem}\label{lem:coladdedvals}
In column $c$, the entries in diagonals $\mathbf{D_{1+2i}}$ and $\mathbf{D_{2+2i}}$ are given by the following pairs for each $1 \le i \le t$:
\begin{equation*}
  n+4ni+2c-1 \quad \text{and} \quad -3n-4ni-2c+1.
\end{equation*}
Similarly, for diagonals $\mathbf{D_{h+2i}}$ and $\mathbf{D_{h+1+2i}}$, the entries are given by the following pairs:
\begin{equation*}
  -n-4ni-2-2((c-h-1) \bmod{n}) \quad \text{and} \quad 3n+4ni+2+2((c-h-1) \bmod{n}).
\end{equation*}
\end{lem}

In preparation for calculating the complete list of column partial sums $S_C(c, p)$, we first determine the list of partial sums generated by the diagonals, $\mathbf{D}_{1+2i}$ and $\mathbf{D}_{2+2i}$. Because the column partial sum calculations begin with this specific pair of diagonals, this initial sequence accounts for the first $2t$ partial sums.

\begin{lem}\label{lem:coladdedsums_D12}
For any column $1 \le c \le n$, summing the entries of the diagonals $\mathbf{D}_{1+2i}$ and $\mathbf{D}_{2+2i}$ top-to-bottom starting in cell $(c+2, c)$ with $1 \le i \le t$ yields the first $2t$ column partial sums $S_C(c, 1), \dots, S_C(c, 2t)$ as the sequence:
\begin{equation*}
5n+2c-1, \ -2n, \ 7n+2c-1, \ -4n, \ \dots, \ 3n+2nt+2c-1, \ -2nt.
\end{equation*}
\end{lem}

We calculate these partial sums by adding the entries in column $c$ defined in \Cref{lem:coladdedvals}. Adding an entry of $\mathbf{D}_{1+2i}$ yields the partial sum $3n+2ni+2c-1$. Adding an entry of $\mathbf{D}_{2+2i}$ brings the cumulative sum to $-2ni$. Expanding this sequence from $i=1$ to $t$ directly results in the $2t$ partial sums above.

Similarly, we determine the partial sums generated by the diagonals, $\mathbf{D}_{h+2i}$ and $\mathbf{D}_{h+1+2i}$. As these sums occur in the interior of the sequence of partial sums, they will be shifted by the preceding partial sum denoted $S_C(c,x) = \mathcal{S}$, to complete the partial sum evaluations for the remaining entries in each column.

\begin{lem}\label{lem:coladdedsums_Dh}
If the current partial sum $S_C(c,x) = \mathcal{S}$, then for any column $1 \le c \le n$, summing the entries of the diagonals $\mathbf{D}_{h+2i}$ and $\mathbf{D}_{h+1+2i}$ top-to-bottom starting in cell $(c+h+1, c)$ with $1 \le i \le t$ yields a sequence of $2t$ partial sums $S_C(c, x+1), \dots, S_C(c, x+2t)$. 

For columns $1 \le c \le h$, the sequence is:
\begin{equation*}
\mathcal{S}-6n-2c+1, \ \mathcal{S}+2n, \ \mathcal{S}-8n-2c+1, \ \mathcal{S}+4n, \ \dots, \ \mathcal{S}-4n-2nt-2c+1, \ \mathcal{S}+2nt.
\end{equation*}

For columns $h+1 \le c \le n$, the sequence is:
\begin{equation*}
\mathcal{S}-4n-2c+1, \ \mathcal{S}+2n, \ \mathcal{S}-6n-2c+1, \ \mathcal{S}+4n, \ \dots, \ \mathcal{S}-2n-2nt-2c+1, \ \mathcal{S}+2nt.
\end{equation*}
\end{lem}

To evaluate these sequences, we take the preceding cumulative sum $\mathcal{S}$ and sequentially add the entries defined in \Cref{lem:coladdedvals}. The calculation splits into two cases based on the column $c$, generating the partial sums shown above.

With these lemmas, we now construct the complete partial sum sequence $S_C(c,p)$. We evaluate the first $k-3$ terms of $S_C(c,p)$ by sequentially adding: the $2t$ entries from $\mathbf{D_{1+2i}}$ and $\mathbf{D_{2+2i}}$ in \Cref{lem:coladdedsums_D12}, the two consecutive base $H(n;5)$ entries (\cref{tab:Hn5colentries}), and the $2t$ entries from $\mathbf{D_{h+2i}}$ and $\mathbf{D_{h+1+2i}}$ in \Cref{lem:coladdedsums_Dh}. Because of the structure of the base $H(n;5)$ array, we must partition these calculations into five cases depending on column $c$: $c=1$, $2 \le c \le h-1$, $c=h$, $h+1 \le c \le n-1$, and $c=n$.

We observe that the partial sums within a given column partition into disjoint intervals: a set of partial sums bounded above by $3n+2nt+2c-1$, and a second set of partial sums close to $2nk$. A direct comparison between the upper bound of the first set and the lower bound of the second set confirms these intervals never overlap. Since the individual sequences comprising these sets are strictly monotonic with respect to $i$ and occupy mutually disjoint sub-intervals, all $k-3$ partial sums are guaranteed to be distinct. As an aid in these comparisons, note that since $k \ge 13$, we have $t=(k-5)/4 \ge 2$. Additionally, $k=4t+5$ implies that $2nk = 8nt+10n$.

By evaluating these cumulative additions and reducing negative sums modulo $2nk+1$, we obtain the following sequences of $k-3$ partial sums for each of the five column cases: 

\paragraph{Column $1$:}
After the first $2t$ partial sums from diagonals $\mathbf{D_{1+2i}}$ and $\mathbf{D_{2+2i}}$, the sum is $S_C(1, 2t) \equiv 2nk-2nt+1$. We add the two base $H(n;5)$ entries for column $1$, which are $-2n+1$ and $n+1$, yielding $S_C(1, 2t+1) = 2nk-2nt-2n+2$ and $S_C(1, 2t+2) = 2nk-2nt-n+3$. Finally, we list the partial sums from diagonals $\mathbf{D}_{h+2i}$ and $\mathbf{D}_{h+1+2i}$, found in \Cref{lem:coladdedsums_Dh} where $\mathcal{S} = S_C(1, 2t+2)$:
\small
\begin{align*}
    & 5n+1, \ 2nk-2n+1, \ 7n+1, \ 2nk-4n+1, \ \dots, \ 3n+2nt+1, \ 2nk-2nt+1, \ \textbf{2nk-2nt-2n+2},\\
    & \textbf{2nk-2nt-n+3}, \ 2nk-2nt-7n+2, \ 2nk-2nt+n+3, \ \dots, \ 2nk-4nt-5n+2, \ 2nk-n+3
\end{align*}
\normalsize

\paragraph{Columns $2 \le c \le h-1$:}
After the first $2t$ partial sums, the sum is $S_C(c, 2t) \equiv 2nk-2nt+1$. We add the two base $H(n;5)$ entries for this column range, which are $-n-2c+2$ and $n+2c-1$, yielding $S_C(c, 2t+1) = 2nk-2nt-n-2c+3$ and $S_C(c, 2t+2) = 2nk-2nt+2$. Finally, we list the partial sums from diagonals $\mathbf{D}_{h+2i}$ and $\mathbf{D}_{h+1+2i}$, found in \Cref{lem:coladdedsums_Dh} where $\mathcal{S} = S_C(c, 2t+2)$:
\small
\begin{align*}
    & 5n+2c-1, \ 2nk-2n+1, \ \dots, \ 3n+2nt+2c-1, \ 2nk-2nt+1, \ \textbf{2nk-2nt-n-2c+3}, \\
    & \textbf{2nk-2nt+2}, \ 2nk-2nt-6n-2c+3, \ 2nk-2nt+2n+2, \ \dots, \ 2nk-4nt-4n-2c+3, \ 1
\end{align*}
\normalsize

\paragraph{Column $h$:}
After the first $2t$ partial sums, the sum is $S_C(h, 2t) \equiv 2nk-2nt+1$. We add the two base $H(n;5)$ entries for column $h$, which are $-2n$ and $2n+2$, yielding $S_C(h, 2t+1) = 2nk-2nt-2n+1$ and $S_C(h, 2t+2) = 2nk-2nt+3$. Finally, we list the partial sums from diagonals $\mathbf{D}_{h+2i}$ and $\mathbf{D}_{h+1+2i}$, found in \Cref{lem:coladdedsums_Dh} where $\mathcal{S} = S_C(h, 2t+2)$:
\small
\begin{align*}
    & 5n+2h-1, \ 2nk-2n+1, \ \dots, \ 3n+2nt+2h-1, \ 2nk-2nt+1, \ \textbf{2nk-2nt-2n+1},\\
    & \textbf{2nk-2nt+3}, \ 2nk-2nt-6n-2h+4, \ 2nk-2nt+2n+3, \ \dots, \ 2nk-4nt-4n-2h+4, \ 2
\end{align*}
\normalsize

\paragraph{Columns $h+1 \le c \le n-1$:}
After the first $2t$ partial sums, the sum is $S_C(c, 2t) \equiv 2nk-2nt+1$. We add the two base $H(n;5)$ entries for this column range, which are $-4n+2c-3$ and $4n-2c+2$, yielding $S_C(c, 2t+1) = 2nk-2nt-4n+2c-2$ and $S_C(c, 2t+2) = 2nk-2nt$. Finally, we list the partial sums from diagonals $\mathbf{D}_{h+2i}$ and $\mathbf{D}_{h+1+2i}$, found in \Cref{lem:coladdedsums_Dh} where $\mathcal{S} = S_C(c, 2t+2)$. Because $c > h$, we use the second sequence of partial sums from \Cref{lem:coladdedsums_Dh}.
\small
\begin{align*}
    & 5n+2c-1, \ 2nk-2n+1, \ \dots, \ 3n+2nt+2c-1, \ 2nk-2nt+1, \ \textbf{2nk-2nt-4n+2c-2},\\
    & \textbf{2nk-2nt}, \ 2nk-2nt-4n-2c+1, \ 2nk-2nt+2n, \ \dots, \ 2nk-4nt-2n-2c+1, \ 2nk
\end{align*}
\normalsize

\paragraph{Column $n$:}
After the first $2t$ partial sums, the sum is $S_C(n, 2t) \equiv 2nk-2nt+1$. We add the two base $H(n;5)$ entries for column $n$, which are $-2n-3$ and $2n+1$, yielding $S_C(n, 2t+1) = 2nk-2nt-2n-2$ and $S_C(n, 2t+2) = 2nk-2nt-1$. Finally, we list the partial sums from diagonals $\mathbf{D}_{h+2i}$ and $\mathbf{D}_{h+1+2i}$, found in \Cref{lem:coladdedsums_Dh} where $\mathcal{S} = S_C(n, 2t+2)$. As in the previous case, we use the second sequence from \Cref{lem:coladdedsums_Dh}:
\small
\begin{align*}
    & 7n-1, \ 2nk-2n+1, \ \dots, \ 5n+2nt-1, \ 2nk-2nt+1, \ \textbf{2nk-2nt-2n-2},\\
    & \textbf{2nk-2nt-1}, \  2nk-2nt-6n, \ 2nk-2nt+2n-1, \ \dots, \ 2nk-4nt-4n, \ 2nk-1
\end{align*}
\normalsize

We note three properties of these calculations. First, the $k-3$ partial sums in any column are distinct and nonzero modulo $2nk+1$. Second, there are only five possible values for the partial sum $S_C(c,k-3)$ in any column: $1$, $2$, $2nk-n+3$, $2nk-1$, and $2nk$. Third, no partial sum from $S_C(c,1)$ to $S_C(c,k-3)$ falls between $3$ and $5n$. 

Using the five possible sums $S_C(c,k-3)$ noted above, we calculate the final three partial sums using the values in \Cref{tab:Hn5colentries}. The results are in \Cref{tab:finalthreesums}. We see that for any column, the final three partial sums are either $0$ or strictly between $3$ and $5n$. Therefore, the columns of Heffter arrays constructed in this manner are simple in their natural ordering. 

\begin{table}[h]
\centering
\caption{Final Three Column Partial Sums}
\label{tab:finalthreesums}
\setlength{\tabcolsep}{4pt} 
\footnotesize
\begin{tabular}{lccccc}
\toprule
 \textbf{Column} & \textbf{$\mathbf{S_C(c, k-3)}$} && $\mathbf{S_C(c, k-2)}$ & $\mathbf{S_C(c, k-1)}$ & $\mathbf{S_C(c, k)}$ \\
\midrule
Column 1 & $2nk-n+3$ && $2n+4$ & $3n+3$ & $0$ \\
\addlinespace
Column 2 & $1$ && $5n-1$ & $4n+1$ & $0$ \\
\midrule
Odd $3 \le c \le h-1$ & $1$ && $4n-(c-1)/2+2$ & $3n+(c+1)/2+2$ & $0$ \\
\addlinespace
Even $4 \le c \le h-2$ & $1$ && $5n-c/2$ &  $4n+c/2$ & $0$ \\
\midrule
Column $h$ & $2$ && $(19n+3)/4$ & $(15n+3)/4$ & $0$ \\
\midrule
Odd $h+1 \le c \le n-2$ & $2nk$ && $4n+(c-1)/2-1$ & $5n-(c+1)/2-1$ & $0$ \\
\addlinespace
Even $h+2 \le c \le n-3$ & $2nk$ && $3n+c/2+1$ & $4n-c/2+1$ & $0$ \\
\midrule
Column $n-1$ & $2nk$ && $(7n+1)/2$ & $3n+1$ & $0$ \\
\addlinespace
Column $n$ & $2nk-1$ && $5n-2$ & $5n-1$ & $0$ \\
\bottomrule
\end{tabular}
\end{table}

\subsection{Rows are simple}

In this section, we prove that the rows of the constructed $H(n;k)$ are simple. 

To prove a row is simple, we compute its partial sums and ensure they are pairwise distinct. For each row $r \in \{1,2,\dots,n\}$, the row partial sums $S_R(r, p)$ are calculated left-to-right, starting in cell $(r, r+3)$ and wrapping around modulo $n$. As before, this choice of starting cell does not affect simplicity.

\subsubsection{Example: Globally Simple $H(19;17)$ Row Partial Sums}
We illustrate this process using the globally simple $H(19;17)$ in \Cref{tab:H1917}. \Cref{tab:H1917_Row_Sums} displays the partial sums for each row. For instance, the values in the first row of the table correspond to the sequence $S_R(1, 1), S_R(1, 2), \dots, S_R(1, 17)$. Shaded cells indicate partial sums derived from the base array $H(19;5)$.

\begin{table}[h]
    \centering
    \caption{Row Partial Sums Modulo $2nk+1=647$ for Globally Simple $H(19;17)$}
    \label{tab:H1917_Row_Sums}
    \setlength{\tabcolsep}{4pt} 
    \scriptsize
        
    \begin{tabular}{lc*{17}{c}} 
    \toprule
    Row 1 && 311 & 36 & 275 & 72 & 239 & 108 & \cellcolor{lightgray}148 & \cellcolor{lightgray}91 & $-$219 $\equiv$ 428 & 55 & 464 & 19 & 500 & 630 & \cellcolor{lightgray}536 & \cellcolor{lightgray}554 & \cellcolor{lightgray}0 \\ 
    Row 2 && 313 & 36 & 277 & 72 & 241 & 108 & \cellcolor{lightgray}164 & \cellcolor{lightgray}109 & 444 & 73 & 480 & 37 & 516 & 1 & \cellcolor{lightgray}588 & \cellcolor{lightgray}571 & \cellcolor{lightgray}0 \\ 
    Row 3 && 315 & 36 & 279 & 72 & 243 & 146 & \cellcolor{lightgray}200 & \cellcolor{lightgray}147 & 480 & 111 & 516 & 75 & 552 & 1 & \cellcolor{lightgray}571 & \cellcolor{lightgray}555 & \cellcolor{lightgray}0 \\ 
    Row 4 && 317 & 36 & 281 & 72 & 207 & 108 & \cellcolor{lightgray}160 & \cellcolor{lightgray}109 & 440 & 73 & 476 & 37 & 550 & 1 & \cellcolor{lightgray}587 & \cellcolor{lightgray}572 & \cellcolor{lightgray}0 \\ 
    Row 5 && 319 & 36 & 283 & 110 & 247 & 146 & \cellcolor{lightgray}196 & \cellcolor{lightgray}147 & 476 & 111 & 512 & 37 & 548 & 1 & \cellcolor{lightgray}570 & \cellcolor{lightgray}556 & \cellcolor{lightgray}0 \\ 
    Row 6 && 321 & 36 & 247 & 72 & 211 & 108 & \cellcolor{lightgray}156 & \cellcolor{lightgray}109 & 436 & 73 & 510 & 37 & 546 & 1 & \cellcolor{lightgray}586 & \cellcolor{lightgray}573 & \cellcolor{lightgray}0 \\ 
    Row 7 && 323 & 74 & 287 & 110 & 251 & 146 & \cellcolor{lightgray}192 & \cellcolor{lightgray}147 & 472 & 73 & 508 & 37 & 544 & 1 & \cellcolor{lightgray}569 & \cellcolor{lightgray}557 & \cellcolor{lightgray}0 \\ 
    Row 8 && 287 & 36 & 251 & 72 & 215 & 108 & \cellcolor{lightgray}152 & \cellcolor{lightgray}109 & 470 & 73 & 506 & 37 & 542 & 1 & \cellcolor{lightgray}585 & \cellcolor{lightgray}574 & \cellcolor{lightgray}0 \\ 
    Row 9 && 289 & 36 & 253 & 72 & 217 & 108 & \cellcolor{lightgray}150 & \cellcolor{lightgray}109 & 468 & 73 & 504 & 37 & 540 & 1 & \cellcolor{lightgray}568 & \cellcolor{lightgray}558 & \cellcolor{lightgray}0 \\ 
    Row 10 && 291 & 36 & 255 & 72 & 219 & 108 & \cellcolor{lightgray}147 & \cellcolor{lightgray}110 & 467 & 74 & 503 & 38 & 539 & 2 & \cellcolor{lightgray}585 & \cellcolor{lightgray}566 & \cellcolor{lightgray}0 \\ 
    Row 11 && 293 & 36 & 257 & 72 & 221 & 108 & \cellcolor{lightgray}128 & \cellcolor{lightgray}107 & 462 & 71 & 498 & 35 & 534 & 646 & \cellcolor{lightgray}574 & \cellcolor{lightgray}582 & \cellcolor{lightgray}0 \\ 
    Row 12 && 295 & 36 & 259 & 72 & 223 & 108 & \cellcolor{lightgray}130 & \cellcolor{lightgray}107 & 460 & 71 & 496 & 35 & 532 & 646 & \cellcolor{lightgray}558 & \cellcolor{lightgray}565 & \cellcolor{lightgray}0 \\ 
    Row 13 && 297 & 36 & 261 & 72 & 225 & 108 & \cellcolor{lightgray}132 & \cellcolor{lightgray}107 & 458 & 71 & 494 & 35 & 530 & 646 & \cellcolor{lightgray}575 & \cellcolor{lightgray}581 & \cellcolor{lightgray}0 \\ 
    Row 14 && 299 & 36 & 263 & 72 & 227 & 108 & \cellcolor{lightgray}134 & \cellcolor{lightgray}107 & 456 & 71 & 492 & 35 & 528 & 646 & \cellcolor{lightgray}559 & \cellcolor{lightgray}564 & \cellcolor{lightgray}0 \\ 
    Row 15 && 301 & 36 & 265 & 72 & 229 & 108 & \cellcolor{lightgray}136 & \cellcolor{lightgray}107 & 454 & 71 & 490 & 35 & 526 & 646 & \cellcolor{lightgray}576 & \cellcolor{lightgray}580 & \cellcolor{lightgray}0 \\ 
    Row 16 && 303 & 36 & 267 & 72 & 231 & 108 & \cellcolor{lightgray}138 & \cellcolor{lightgray}107 & 452 & 71 & 488 & 35 & 524 & 646 & \cellcolor{lightgray}560 & \cellcolor{lightgray}563 & \cellcolor{lightgray}0 \\ 
    Row 17 && 305 & 36 & 269 & 72 & 233 & 108 & \cellcolor{lightgray}140 & \cellcolor{lightgray}107 & 450 & 71 & 486 & 35 & 522 & 646 & \cellcolor{lightgray}577 & \cellcolor{lightgray}579 & \cellcolor{lightgray}0 \\ 
    Row 18 && 307 & 36 & 271 & 72 & 235 & 108 & \cellcolor{lightgray}142 & \cellcolor{lightgray}107 & 448 & 71 & 484 & 35 & 520 & 646 & \cellcolor{lightgray}561 & \cellcolor{lightgray}552 & \cellcolor{lightgray}0 \\ 
    Row 19 && 309 & 36 & 273 & 72 & 237 & 108 & \cellcolor{lightgray}144 & \cellcolor{lightgray}106 & 445 & 70 & 481 & 34 & 517 & 645 & \cellcolor{lightgray}587 & \cellcolor{lightgray}588 & \cellcolor{lightgray}0 \\ 
    \bottomrule
    \end{tabular}

\end{table}

\subsubsection{Row Partial Sums Are Distinct Modulo $2nk+1$}

As with the columns, we state several lemmas used to derive the row partial sums for specific sets of rows.

\begin{lem}
Let $H(n;5)$ be the base array generated via the Dinitz and Wanless construction \cite{dinitz2017}. The non-empty entries of this array are exactly those detailed in \Cref{tab:Hn5rowentries}, positioned such that in each row $r$, two consecutive entries start at cell $(r, r+h-1)$ and three consecutive entries start at cell $(r, r-1)$.
\end{lem}

\begin{table}[h]
\centering
\caption{Base $H(n;5)$ Row Entries}
\label{tab:Hn5rowentries}
\setlength{\tabcolsep}{4pt} 
\footnotesize
\begin{tabular}{lcccccc}
\toprule
\textbf{Row} & \multicolumn{2}{c}{\textbf{Two Consecutive Entries}} && \multicolumn{3}{c}{\textbf{Three Consecutive Entries}} \\
\midrule
Row 1 & $2n+2$, & $-3n$, && $-5n+1$, & $n-1$, & $5n-2$ \\
\addlinespace
Row 2 & $3n-1$, & $-3n+2$, && $-3n-3$, & $-n+2$, & $4n$   \\
\midrule
Odd $3 \le r \le h-1$ & $3n-2r+3$, & $-3n+2r-2$, && $-4n-(r-1)/2$, & $-n+r$, & $5n-(r-1)/2-2$ \\
\addlinespace
Even $4 \le r \le h-2$ & $3n-2r+3$, & $-3n+2r-2$, && $-3n-r/2-2$, & $-n+r$, & $4n-r/2+1$ \\
\midrule
Row $h$ & $2n+1$, & $-2n+1$, && $-(13n+9)/4$, & $-n$, & $(17n+1)/4$ \\
\midrule
Odd $h+1 \le r \le n-2$ & $2r-2$, & $-2r+1$, && $-4n+(r+1)/2-2$, & $n-r$, & $3n+(r+1)/2+2$ \\
\addlinespace
Even $h+2 \le r \le n-3$ & $2r-2$, & $-2r+1$, && $-5n+r/2+1$, & $n-r$, & $4n+r/2$ \\
\midrule
Row $n-1$ & $2n-4$, & $-2n+3$, && $(-9n+1)/2$, & $-(n-1)/2$, & $5n$ \\
\addlinespace
Row $n$ & $2n-2$, & $-2n$, && $-3n-1$, & $1$, & $3n+2$ \\
\bottomrule
\end{tabular}
\end{table}

\begin{lem}
The entries in diagonals $\mathbf{D_{h+1+2i}}$ and $\mathbf{D_{h+2i}}$ for a fixed row $r$ are given by the pairs with $1 \le i \le t$:
\begin{equation*}
  3n+4ni+2+2((r-2-2i) \bmod{n}) \quad \text{and} \quad -n-4ni-2-2((r-1-2i) \bmod{n}).
\end{equation*}
Similarly, the entries in row $r$ for diagonals $\mathbf{D_{2+2i}}$ and $\mathbf{D_{1+2i}}$ are given by the pairs:
\begin{equation*}
  -3n-4ni-1-2((r-2-2i) \bmod{n}) \quad \text{and} \quad n+4ni+1+2((r-1-2i) \bmod{n}).
\end{equation*}
\end{lem}

\begin{lem}\label{lem:rowaddedsums_Dh}
For any row $1 \le r \le n$, summing the entries in the diagonals $\mathbf{D_{h+1+2i}}$ and $\mathbf{D_{h+2i}}$ left-to-right starting in cell $(r, r+3)$ with $1 \le i \le t$ yields the first $2t$ row partial sums $S_R(r, 1), \dots, S_R(r, 2t)$. The pair of partial sums generated at step $i$ is determined by four cases based on the row $r$:
\begin{description}
    \item[Case 1: Rows 1 and 2.]
    The pair of sums is: $7n+4nt-4t+2r-4-(2n-2)i \quad \text{and} \quad (2n-2)i$.
    \item[Case 2: Odd Rows $3 \le r \le (k-3)/2$.]
    For $i < t-(r-3)/2$, the pair is $7n+4nt-4t+2r-4-(2n-2)i$ and $(2n-2)i$. For $i \ge t-(r-3)/2$, the pair is $7n+4nt-4t+2r-4-(2n-2)i$ and $2n+(2n-2)i$.
    \item[Case 3: Even Rows $4 \le r \le (k-5)/2$.]
    For $i \le t-(r-2)/2$, the pair is $7n+4nt-4t+2r-4-(2n-2)i$ and $(2n-2)i$. For $i > t-(r-2)/2$, the pair is $5n+4nt-4t+2r-4-(2n-2)i$ and $(2n-2)i$.
    \item[Case 4: Rows $(k-1)/2 \le r \le n$.]
    The pair of sums is: $5n+4nt-4t+2r-4-(2n-2)i \quad \text{and} \quad (2n-2)i$.
\end{description}
\end{lem}

\begin{lem}\label{lem:rowaddedsums_D12}
For any row $1 \le r \le n$, summing the entries of diagonals $\mathbf{D_{2+2i}}$ and $\mathbf{D_{1+2i}}$ left-to-right starting in cell $(r, r+h+2)$ with $1 \le i \le t$ yields a sequence of $2t$ partial sums $S_R(r, x+1), \dots, S_R(r, x+2t)$, where $\mathcal{S} = S_R(r,x)$ is the preceding partial sum. The pair of partial sums generated at step $i$ is determined by four cases based on the row $r$:
\begin{description}
    \item[Case 1: Rows 1 and 2.]
    The pair of sums is: $\mathcal{S}-7n-4nt+4t-2r+5+(2n-2)i \quad \text{and} \quad \mathcal{S}-(2n-2)i$.

    \item[Case 2: Odd Rows $3 \le r \le (k-3)/2$.]
    For $i < t-(r-3)/2$, the pair is $\mathcal{S}-7n-4nt+4t-2r+5+(2n-2)i$ and $\mathcal{S}-(2n-2)i$. For $i \ge t-(r-3)/2$, the pair is $\mathcal{S}-7n-4nt+4t-2r+5+(2n-2)i$ and $\mathcal{S}-2n-(2n-2)i$.
    
    \item[Case 3: Even Rows $4 \le r \le (k-5)/2$.]
    For $i \le t-(r-2)/2$, the pair is $\mathcal{S}-7n-4nt+4t-2r+5+(2n-2)i$ and $\mathcal{S}-(2n-2)i$. For $i > t-(r-2)/2$, the pair is $\mathcal{S}-5n-4nt+4t-2r+5+(2n-2)i$ and $\mathcal{S}-(2n-2)i$.
    
    \item[Case 4: Rows $(k-1)/2 \le r \le n$.]
    The pair of sums is: $\mathcal{S}-5n-4nt+4t-2r+5+(2n-2)i \quad \text{and} \quad \mathcal{S}-(2n-2)i$.

\end{description}
\end{lem}

Using these lemmas, we construct the complete partial sum sequence $S_R(r,p)$. We evaluate the first $k-3$ terms of $S_R(r,p)$ by sequentially adding: the $2t$ entries from diagonals $\mathbf{D_{h+1+2i}}$ and $\mathbf{D_{h+2i}}$ in \Cref{lem:rowaddedsums_Dh}, the two consecutive base $H(n;5)$ entries (\cref{tab:Hn5rowentries}), and the $2t$ entries from diagonals $\mathbf{D_{2+2i}}$ and $\mathbf{D_{1+2i}}$ in \Cref{lem:rowaddedsums_D12}. The structure of the base $H(n;5)$ array requires us to partition these sums into eight cases depending on row $r$: $r=1$, $2$, odd and even rows from $3 \le r \le (k-3)/2$, $(k-1)/2 \le r \le h-1$, $r=h$, $h+1 \le r \le n-1$, and $r=n$.

To see that the partial sums within each row case are distinct, observe that they partition into disjoint intervals: a set of partial sums bounded above by $5n+2nt-2t-3$, a second set of partial sums bounded between $5n+2nt-2t+4$ and $5n+4nt-4t+2r-2$, and a third set of partial sums close to $2nk$. Since the individual sequences comprising these sets are strictly monotonic with respect to $i$ and occupy mutually disjoint sub-intervals, all $k-3$ partial sums are guaranteed to be distinct. Recall that $k \ge 13$ implies $t=(k-5)/4 \ge 2$, and $k=4t+5$ implies $2nk = 8nt+10n$. 

By evaluating these cumulative additions and reducing negative sums modulo $2nk+1$, we obtain the following sequences of $k-3$ partial sums for each row case:

\paragraph{Row $1$:}
After the first $2t$ partial sums, the sum is $S_R(1, 2t) \equiv 2nt-2t$. We add the two base $H(n;5)$ entries for row $1$, which are $2n+2$ and $-3n$, yielding $S_R(1, 2t+1) = 2nt+2n-2t+2$ and $S_R(1, 2t+2) = 2nt-n-2t+2$. Finally, we list the subsequent partial sums from the remaining diagonals found in \Cref{lem:rowaddedsums_Dh} where for each row $1 \le r \le n$, $\mathcal{S} = S_R(r, 2t+2)$:
\small
\begin{align*}
& 5n+4nt-4t, \ 2n-2, \ 3n+4nt-4t+2, \ 4n-4, \ \dots, \ 7n+2nt-2t-2, \ 2nt-2t, \ \textbf{2nt+2n-2t+2},\\
& \textbf{2nt-n-2t+2}, \ 2nk-2nt-6n+2t+4, \ 2nt-3n-2t+4, \ \dots, \ 2nk-8n+6, \ 2nk-n+3
\end{align*} 
\normalsize

\paragraph{Row $2$:}
After the first $2t$ partial sums, the sum is $S_R(2, 2t) \equiv 2nt-2t$. We add the two base $H(n;5)$ entries for row $2$, which are $3n-1$ and $-3n+2$, yielding $S_R(2, 2t+1) = 2nt+3n-2t-1$ and $S_R(2, 2t+2) = 2nt-2t+1$. Finally, we list the subsequent partial sums:
\small
\begin{align*}
& 5n+4nt-4t+2, \ 2n-2, \ 3n+4nt-4t+4, \ 4n-4, \ \dots, \ 7n+2nt-2t, \ 2nt-2t, \  \textbf{2nt+3n-2t-1},\\
& \textbf{2nt-2t+1}, \ 2nk-2nt-5n+2t+1, \ 2nt-2n-2t+3, \ \dots, \ 2nk-7n+3, \ 1
\end{align*} 
\normalsize 

Rows $3$ through $(k-3)/2$ contain the starting cells of diagonals $\mathbf{D_{h+1+2i}}$, $\mathbf{D_{h+2i}}$, $\mathbf{D_{2+2i}}$, and $\mathbf{D_{1+2i}}$. This causes the sequence of partial sums in these rows to jump when adding the entry in the starting cell. As detailed in \Cref{lem:rowaddedsums_Dh}, when adding the entries of $\mathbf{D_{h+1+2i}}$ and $\mathbf{D_{h+2i}}$ in a given odd row with $3 \le r \le (k-3)/2$, as long as $i < t-(r-3)/2$, after adding an entry from $\mathbf{D_{h+2i}}$, the current partial sum will be $(2n-2)i$. However, when $i = t-(r-3)/2$, the entry we are adding is the starting cell of $\mathbf{D_{h+2i}}$. For this and all subsequent $i$ in this row, adding an entry in $\mathbf{D_{h+2i}}$ will yield the partial sum $2n+(2n-2)i$. When listing the partial sums in these rows, we use set notation, $\{(2n-2)i, 2n+(2n-2)i\}$, to represent the two possible partial sums. Depending on $r$ and $i$, only one of these partial sums will be realized. Because each of the potential sums in these sets is strictly distinct from the rest of the partial sums in the row, it is not necessary to split this into cases. The potential partial sums for $\mathbf{D_{2+2i}}$ and $\mathbf{D_{1+2i}}$ are detailed in \Cref{lem:rowaddedsums_D12}.

\paragraph{Odd Rows $3 \le r \le (k-3)/2$:}
After the first $2t$ partial sums, the sum is $S_R(r, 2t) \equiv 2n+2nt-2t$. We add the two base $H(n;5)$ entries for this row range, which are $3n-2r+3$ and $-3n+2r-2$, yielding $S_R(r, 2t+1) = 2nt+5n-2t-2r+3$ and $S_R(r, 2t+2) = 2nt+2n-2t+1$. For partial sums that vary depending on $i$ and the row, as noted in \Cref{lem:rowaddedsums_D12}, we use set notation to list the possible values. Finally, we list the subsequent partial sums:
\small
\begin{align*}
& 5n+4nt-4t+2r-2, \ \{2n-2, 4n-2\}, \ 3n+4nt-4t+2r, \ \{4n-4, 6n-4\}, \ \dots \\
& \dots, \ 7n+2nt-2t+2r-4, \ 2n+2nt-2t, \ \textbf{2nt+5n-2t-2r+3}, \\
& \textbf{2nt+2n-2t+1}, \  2nk-2nt-3n+2t-2r+5, \ \{2nt-2n-2t+3, 2nt-2t+3\} \dots, \ 2nk-5n-2r+7, \ 1
\end{align*}
\normalsize

\paragraph{Even Rows $4 \le r \le (k-5)/2$:}
After the first $2t$ partial sums, the sum is $S_R(r, 2t) \equiv 2nt-2t$. We add the two base $H(n;5)$ entries for this row range, which are $3n-2r+3$ and $-3n+2r-2$, yielding $S_R(r, 2t+1) = 2nt+3n-2t-2r+3$ and $S_R(r, 2t+2) = 2nt-2t+1$. As with the odd rows, we use set notation to represent the partial sums that vary depending on $i$ and the row, which sufficiently demonstrates the distinctness of the sequence of partial sums:
\small
\begin{align*}
& \{5n+4nt-4t+2r-2, 7n+4nt-4t+2r-2\}, \ 2n-2, \dots \{5n+2nt-2t+2r-4, 7n+2nt-2t+2r-4\}, \\
& \dots, 2nt-2t, \ \textbf{2nt+3n-2t-2r+3}, \ \textbf{2nt-2t+1}, \\
& \{2nk-2nt-5n+2t-2r+5, 2nk-2nt-7n+2t-2r+5\}, \ 2nt-2n-2t+3, \dots, 2nk-5n-2r+7, \ 1
\end{align*}
\normalsize

\paragraph{Rows $(k-1)/2 \le r \le h-1$:}
After the first $2t$ partial sums, the sum is $S_R(r, 2t) \equiv 2nt-2t$. We add the two base $H(n;5)$ entries for this row range, which are $3n-2r+3$ and $-3n+2r-2$, yielding $S_R(r, 2t+1) = 3n+2nt-2t-2r+3$ and $S_R(r, 2t+2) = 2nt-2t+1$. Finally, we list the subsequent partial sums:
\small
\begin{align*}
& 3n+4nt-4t+2r-2, \ 2n-2, \ n+4nt-4t+2r, \ 4n-4, \dots, 5n+2nt-2t+2r-4, \ 2nt-2t, \\
& \textbf{3n+2nt-2t-2r+3}, \ \textbf{2nt-2t+1}, \ 2nk-2nt-3n+2t-2r+5, \ 2nt-2n-2t+3, \dots, 2nk-5n-2r+7, \ 1
\end{align*}
\normalsize

\paragraph{Row $h$:}
After the first $2t$ partial sums, the sum is $S_R(h, 2t) \equiv 2nt-2t$. We add the two base entries for row $h$, which are $2n+1$ and $-2n+1$, yielding $S_R(h, 2t+1) = 2n+2nt-2t+1$ and $S_R(h, 2t+2) = 2nt-2t+2$. Finally, we list the subsequent partial sums:
\small
\begin{align*}
& 3n+4nt-4t+2h-2, \ 2n-2, \ n+4nt-4t+2h, \ 4n-4, \dots, 5n+2nt-2t+2h-4, \ 2nt-2t, \\
& \textbf{2n+2nt-2t+1}, \ \textbf{2nt-2t+2}, \ 2nk-2nt-3n+2t-2h+6, \ 2nt-2n-2t+4, \dots, 2nk-5n-2h+8, \ 2
\end{align*}
\normalsize

\paragraph{Rows $h+1 \le r \le n-1$:}
After the first $2t$ partial sums, the sum is $S_R(r, 2t) \equiv 2nt-2t$. We add the two base entries for this row range, which are $2r-2$ and $-2r+1$, yielding $S_R(r, 2t+1) = 2nt-2t+2r-2$ and $S_R(r, 2t+2) = 2nt-2t-1$. Finally, we list the subsequent partial sums:
\small
\begin{align*}
& 3n+4nt-4t+2r-2, \ 2n-2, \ n+4nt-4t+2r, \ 4n-4, \dots, 5n+2nt-2t+2r-4, \ 2nt-2t, \\
& \textbf{2nt-2t+2r-2}, \ \textbf{2nt-2t-1}, \ 2nk-2nt-3n+2t-2r+3, \ 2nt-2n-2t+1, \dots, 2nk-5n-2r+5, \ 2nk 
\end{align*}
\normalsize

\paragraph{Row $n$:}
After the first $2t$ partial sums, the sum is $S_R(n, 2t) \equiv 2nt-2t$. We add the two base entries for row $n$, which are $2n-2$ and $-2n$, yielding $S_R(n, 2t+1) = 2n+2nt-2t-2$ and $S_R(n, 2t+2) = 2nt-2t-2$. Finally, we list the subsequent partial sums:
\begin{align*}
& 5n+4nt-4t-2, \ 2n-2, \ 3n+4nt-4t, \ 4n-4, \dots, 7n+2nt-2t-4, \ 2nt-2t, \\
& \textbf{2n+2nt-2t-2}, \ \textbf{2nt-2t-2}, \ 2nk-2nt-5n+2t+2, \ 2nt-2n-2t, \dots, 2nk-7n+4, \ 2nk-1 
\end{align*}

We note three properties of these calculations that allow us to easily verify the distinctness of the partial sums. First, the initial $k-3$ partial sums in any row are distinct and nonzero modulo $2nk+1$. Second, there are only five possible values for the partial sum $S_R(r, k-3)$ in any row: $1$, $2$, $2nk$, $2nk-1$, and $2nk-n+3$. Third, for any row $r$, no partial sum $S_R(r, p)$ (for $1 \le p \le k-3$) falls between $2nk-5n-2r+8$ and $2nk-3n$.

Using these properties, we calculate the final three partial sums using the values in \Cref{tab:Hn5rowentries}, yielding the results compiled in \Cref{tab:finalthreerowsums}. We observe that for the general rows $2 \le r \le n-2$ and row $n$, these final three partial sums, $S_R(r, k-2), S_R(r, k-1)$, and $S_R(r, k)$, are either $0$ or fall strictly within the interval $[2nk-5n-2r+8, 2nk-3n]$. Because this interval is completely disjoint from the values attained by the first $k-3$ partial sums, no clashing can occur between the two sections of these rows.

For rows $1$ and $n-1$, distinctness is similarly verified by observing the largest values among the first $k-3$ partial sums. In row $1$ the largest two values among the first $k-3$ partial sums are $2nk-8n+6$ and $2nk-n+3$. The final three partial sums are either $0$ or fall strictly between these two values. In row $n-1$ the largest two values from the first $k-3$ partial sums are $2nk-7n+7$ and $2nk$. The final three partial sums are either $0$ or fall strictly between these two values. Because these final sums are bounded between the largest and second-largest preceding sums, they are distinct. Therefore, all partial sums within every row are distinct. This, taken together with the fact that the columns are simple proves that the Heffter arrays constructed in this manner are globally simple.

\begin{table}[h]
\centering
\caption{Final Three Row Partial Sums}
\label{tab:finalthreerowsums}
\setlength{\tabcolsep}{4pt} 
\footnotesize
\begin{tabular}{lccccc}
\toprule
\textbf{Row} & \textbf{$\mathbf{S_R(r, k-3)}$} && $\mathbf{S_R(r, k-2)}$ & $\mathbf{S_R(r, k-1)}$ & $\mathbf{S_R(r, k)}$ \\
\midrule
Row 1 & $2nk-n+3$ && $2nk-6n+4$ & $2nk-5n+3$ & $0$ \\
\addlinespace
Row 2 & $1$ && $2nk-3n-1$ & $2nk-4n+1$ & $0$ \\
\midrule
Odd $3 \le r \le h-1$ & $1$ && $2nk-4n-(r-1)/2+2$ & $2nk-5n+(r+1)/2+2$ & $0$ \\
\addlinespace
Even $4 \le r \le h-2$ & $1$ && $2nk-3n-r/2$ & $2nk-4n+r/2$ & $0$ \\
\midrule
Row $h$ & $2$ && $2nk-(13n+9)/4+3$ & $2nk-(17n+9)/4+3$ & $0$ \\
\midrule
Odd $h+1 \le r \le n-2$ & $2nk$ && $2nk-4n+(r+1)/2-2$ & $2nk-3n-(r-1)/2-2$ & $0$ \\
\addlinespace
Even $h+2 \le r \le n-3$ & $2nk$ && $2nk-5n+r/2+1$ & $2nk-4n-r/2+1$ & $0$ \\
\midrule
Row $n-1$ & $2nk$ && $2nk-(9n-1)/2$ & $2nk-5n+1$ & $0$ \\
\addlinespace
Row $n$ & $2nk-1$ && $2nk-3n-2$ & $2nk-3n-1$ & $0$ \\
\bottomrule
\end{tabular}
\end{table}

\section{$H(n;k)$ Construction for $n \equiv 0 \pmod{4}$ and $k \equiv 1 \pmod{4}$} \label{sec:construction01}

We construct globally simple Heffter arrays $H(n;k)$ for all $n \equiv 0 \pmod{4}$ and $k \equiv 1 \pmod{4}$ with $n > k \geq 13$. The construction proceeds by initializing a base $H(n;5)$ array, as described by Dinitz and Wanless \cite{dinitz2017}, and augmenting with sets of 4 diagonals until $k$ is reached.

\subsection{The Base Array $H(n;5)$ when $n \equiv 0 \pmod{4}$}

Let $h=n/2$ and $q=n/4$. The base $H(n;5)$ is constructed using the following diagonals and ad hoc entries.

\paragraph{Diagonal Assignments}
\begin{align*}
    &\text{A: diag}(h+1, h+2, (h-2), 1, -1, h-3) & &\text{H: diag}(h+2,h+2,-(5n-q-4),2,1,q-1) \\
    &\text{B: diag}(1,2,-(n-2),1,1,h-1) & &\text{I: diag}(h,h+2,(4n+q-1),2,1,q-1)\\
    &\text{C: diag}(1,3,(4n-2),2,-1,q) & &\text{J: diag}(h+1,h+3,(3n+q+1),2,1,q-1)\\
    &\text{D: diag}(2,2,-(4n-1),2,-1,q) & &\text{K: diag}(h+1,1,(n-1),1,2,h-1)\\
    &\text{E: diag}(2,4,(5n-5),2,-1,q-1) & &\text{L: diag}(h+1,2,-(n),1,-2,h-2)\\
    &\text{F: diag}(3,3,-(3n+2),2,-1,q-1) & &\text{M: diag}(2,h+2,(3n-4),1,-2,h-2)\\
    &\text{G: diag}(h+1,h+1,-(4n-q-2),2,1,q-1) & &\text{N: diag}(1,h+2,-(3n-3),1,2,h-1)
\end{align*}

\paragraph{Ad Hoc Cell Values}
The following specific entries are set directly:

\begin{align*}
    &H[1,1] = 3n-2, &&H[1,h+1] = -3n-1, \\
    &H[h,1] = -2n+2, &&H[h,h+1] = 5n-3, &&H[h,n] = -3n, \\
    &H[n-2,n-1] = -h+1, &&H[n-2,n] = 5n-2, \\
    &H[n-1,1] = -2n, &&H[n-1,h] = 5n, \\
    &H[n-1,n-1]=-3n+1, &&H[n-1,n]=-2n+4, \\
    &H[n,1] = 1, &&H[n,2] = 5n-4, &&H[n,h] = -5n+1, \\
    &H[n,h+1] = -2n+3, &&H[n,n]=2n-1. \\
\end{align*}

\subsection{Augmentation with $k-5$ New Diagonals}

To increase the number of filled cells in each row and column from 5 to $k$, we add diagonals in groups of 4. Let $t= (k-5)/4$. For each $i \in \{1, 2, \dots, t\}$, we define the following four diagonals:

\begin{equation*}
\begin{array}{lll}
\toprule
\textbf{Diagonal} & \textbf{Definition} & \textbf{Support} \\
\midrule
\mathbf{D_{2+2i}} & \text{diag}(2+2i, 1, (-n - 4ni - 1), 1, -4, n) & \{n+4ni+1, n+4ni+5, \dots, 5n+4ni-3\} \\
\addlinespace
\mathbf{D_{3+2i}} & \text{diag}(3+2i, 1, (n + 4ni + 3 ), 1, 4, n) & \{n+4ni+3, n+4ni+7, \dots, 5n+4ni-1\} \\
\addlinespace
\mathbf{D_{h+1+2i}} & \text{diag}(2+2i, h+3, (n + 4ni + 2), 1, 4, n) & \{n+4ni+2, n+4ni+6, \dots, 5n+4ni-2\} \\
\addlinespace
\mathbf{D_{h+2+2i}} & \text{diag}(3+2i, h+3, (-n - 4ni - 4), 1, -4, n) & \{n+4ni+4, n+4ni+8, \dots, 5n+4ni\} \\
\bottomrule
\end{array}
\end{equation*}

Combined with the base support $\{1, \dots, 5n\}$ from the $H(n;5)$, the constructed $H(n;k)$ has the required support $\{1, 2, \dots, nk\}$.

\subsection{Row and Column Sums}

The rows and columns of the base $H(n;5)$ sum to zero. We verify the above diagonals do as well. Consider the contribution of the $i$-th set of diagonals to any fixed row $r$. The row $r$ intersects $\mathbf{D_{2+2i}}$ at some column $c$ with value $v$, and intersects $\mathbf{D_{h+1+2i}}$ at column $c'$ with value $v'$. By the diagonal definitions:
\[ v + v' = (-n-4ni-1-4((9r-2-2i) \bmod{n})) + (n+4ni+2+4((r-2-2i) \bmod{n})) = 1. \]
Similarly, the values from $\mathbf{D_{3+2i}}$ and $\mathbf{D_{h+2+2i}}$ sum to $-1$. Thus, the total contribution of these diagonals to row $r$ is $1 - 1 = 0$. For any fixed column $c$, the values from $\mathbf{D_{2+2i}}$ and $\mathbf{D_{3+2i}}$ sum to $2$, while the values from $\mathbf{D_{h+1+2i}}$ and $\mathbf{D_{h+2+2i}}$ sum to $-2$. Thus, the total contribution of these diagonals to column $c$ is $2 - 2 = 0$.

Since the support is correct and every row and column sum is zero, this construction produces a valid integer Heffter array $H(n;k)$.

\subsection{Example: Globally Simple $H(20;17)$}

The 5 light gray diagonals in \cref{tab:H2017} are the base $H(n;5)$. The starting cells for the 12 additional diagonals are highlighted in dark gray.

\begin{table}[h]
    \centering
    \caption{Heffter Array $H(20;17)$}
    \label{tab:H2017}
    \setlength{\tabcolsep}{4pt}
    \scriptsize
        \begin{tabular}{|*{20}{c|}} 
            \hline
\cellcolor{black!15}58 & \cellcolor{black!15}-18 & \cellcolor{black!15}78 & & -312 & 314 & -240 & 242 & -168 & 170 & \cellcolor{black!15}-61 & \cellcolor{black!15}-57 & 311 & -313 & 239 & -241 & 167 & -169 & & \\ \hline
            & \cellcolor{black!15}-79 & \cellcolor{black!15}-17 & \cellcolor{black!15}95 & & -316 & 318 & -244 & 246 & -172 & 174 & \cellcolor{black!15}56 & \cellcolor{black!15}-55 & 315 & -317 & 243 & -245 & 171 & -173 & \\ \hline
            & & \cellcolor{black!15}-62 & \cellcolor{black!15}-16 & \cellcolor{black!15}77 & & -320 & 322 & -248 & 250 & -176 & 178 & \cellcolor{black!15}54 & \cellcolor{black!15}-53 & 319 & -321 & 247 & -249 & 175 & -177 \\ \hline
            \cellcolor{black!40}-101 & & & \cellcolor{black!15}-80 & \cellcolor{black!15}-15 & \cellcolor{black!15}94 & & -324 & 326 & -252 & 254 & -180 & \cellcolor{black!40}102 & \cellcolor{black!15}52 & \cellcolor{black!15}-51 & 323 & -325 & 251 & -253 & 179 \\ \hline
            \cellcolor{black!40}103 & -105 & & & \cellcolor{black!15}-63 & \cellcolor{black!15}-14 & \cellcolor{black!15}76 & & -328 & 330 & -256 & 258 & \cellcolor{black!40}-104 & 106 & \cellcolor{black!15}50 & \cellcolor{black!15}-49 & 327 & -329 & 255 & -257 \\ \hline
            \cellcolor{black!40}-181 & 107 & -109 & & & \cellcolor{black!15}-81 & \cellcolor{black!15}-13 & \cellcolor{black!15}93 & & -332 & 334 & -260 & \cellcolor{black!40}182 & -108 & 110 & \cellcolor{black!15}48 & \cellcolor{black!15}-47 & 331 & -333 & 259 \\ \hline
            \cellcolor{black!40}183 & -185 & 111 & -113 & & & \cellcolor{black!15}-64 & \cellcolor{black!15}-12 & \cellcolor{black!15}75 & & -336 & 338 &\cellcolor{black!40} -184 & 186 & -112 & 114 & \cellcolor{black!15}46 & \cellcolor{black!15}-45 & 335 & -337 \\ \hline
            \cellcolor{black!40}-261 & 187 & -189 & 115 & -117 & & & \cellcolor{black!15}-82 & \cellcolor{black!15}-11 & \cellcolor{black!15}92 & & -340 &\cellcolor{black!40} 262 & -188 & 190 & -116 & 118 & \cellcolor{black!15}44 & \cellcolor{black!15}-43 & 339 \\ \hline
            \cellcolor{black!40}263 & -265 & 191 & -193 & 119 & -121 & & & \cellcolor{black!15}-65 & \cellcolor{black!15}-10 & \cellcolor{black!15}74 & & \cellcolor{black!40}-264 & 266 & -192 & 194 & -120 & 122 & \cellcolor{black!15}42 & \cellcolor{black!15}-41 \\ \hline
            \cellcolor{black!15}-38 & 267 & -269 & 195 & -197 & 123 & -125 & & & \cellcolor{black!15}-83 & \cellcolor{black!15}97 & \cellcolor{black!15}84 & & -268 & 270 & -196 & 198 & -124 & 126 & \cellcolor{black!15}-60 \\ \hline
            \cellcolor{black!15}19 & \cellcolor{black!15}-20 & 271 & -273 & 199 & -201 & 127 & -129 & & & \cellcolor{black!15}-73 & \cellcolor{black!15}8 & \cellcolor{black!15}66 & & -272 & 274 & -200 & 202 & -128 & 130 \\ \hline
            134 & \cellcolor{black!15}21 & \cellcolor{black!15}-22 & 275 & -277 & 203 & -205 & 131 & -133 & & & \cellcolor{black!15}-91 & \cellcolor{black!15}7 & \cellcolor{black!15}85 & & -276 & 278 & -204 & 206 & -132 \\ \hline
            -136 & 138 & \cellcolor{black!15}23 & \cellcolor{black!15}-24 & 279 & -281 & 207 & -209 & 135 & -137 & & & \cellcolor{black!15}-72 & \cellcolor{black!15}6 & \cellcolor{black!15}67 & & -280 & 282 & -208 & 210 \\ \hline
            214 & -140 & 142 & \cellcolor{black!15}25 & \cellcolor{black!15}-26 & 283 & -285 & 211 & -213 & 139 & -141 & & & \cellcolor{black!15}-90 & \cellcolor{black!15}5 & \cellcolor{black!15}86 & & -284 & 286 & -212 \\ \hline
            -216 & 218 & -144 & 146 & \cellcolor{black!15}27 & \cellcolor{black!15}-28 & 287 & -289 & 215 & -217 & 143 & -145 & & & \cellcolor{black!15}-71 & \cellcolor{black!15}4 & \cellcolor{black!15}68 & & -288 & 290 \\ \hline
            294 & -220 & 222 & -148 & 150 & \cellcolor{black!15}29 & \cellcolor{black!15}-30 & 291 & -293 & 219 & -221 & 147 & -149 & & & \cellcolor{black!15}-89 & \cellcolor{black!15}3 & \cellcolor{black!15}87 & & -292 \\ \hline
            -296 & 298 & -224 & 226 & -152 & 154 & \cellcolor{black!15}31 & \cellcolor{black!15}-32 & 295 & -297 & 223 & -225 & 151 & -153 & & & \cellcolor{black!15}-70 & \cellcolor{black!15}2 & \cellcolor{black!15}69 & \\ \hline
            & -300 & 302 & -228 & 230 & -156 & 158 & \cellcolor{black!15}33 & \cellcolor{black!15}-34 & 299 & -301 & 227 & -229 & 155 & -157 & & & \cellcolor{black!15}-88 & \cellcolor{black!15}-9 & \cellcolor{black!15}98 \\ \hline
            \cellcolor{black!15}-40 & & -304 & 306 & -232 & 234 & -160 & 162 & \cellcolor{black!15}35 & \cellcolor{black!15}100 & 303 & -305 & 231 & -233 & 159 & -161 & & & \cellcolor{black!15}-59 & \cellcolor{black!15}-36 \\ \hline
            \cellcolor{black!15}1 & \cellcolor{black!15}96 & & -308 & 310 & -236 & 238 & -164 & 166 & \cellcolor{black!15}-99 & \cellcolor{black!15}-37 & 307 & -309 & 235 & -237 & 163 & -165 & & & \cellcolor{black!15}39 \\ \hline
\end{tabular}
\end{table}

\section{$H(n;k)$ is globally simple for $n \equiv 0 \pmod{4}$ and $k \equiv 1 \pmod{4}$} \label{sec:proof_simple01}

In this section, we prove \Cref{thm:n0k1} by showing the integer Heffter arrays constructed in \Cref{sec:construction01} are globally simple.

\subsection{Columns are simple}

To prove a column is simple, we compute its partial sums and ensure they are pairwise distinct. For each column $c \in \{1,2,\dots,n\}$, the column partial sums $S_C(c, p)$ are calculated top-to-bottom, starting in cell $(c+3, c)$. As in Section 3, this choice of starting cell does not affect simplicity. 

\subsubsection{Example: Globally Simple $H(20;17)$ Column Partial Sums}
We illustrate this process using the globally simple $H(20;17)$ in \Cref{tab:H2017}. \Cref{tab:H2017colsums} displays the computed partial sums for each column. For instance, the values in the first row of the table correspond to the sequence $S_C(1, 1), S_C(1, 2), \dots, S_C(1, 17)$. Shaded cells indicate partial sums derived from the base array $H(20;5)$.

\begin{table}[h]
    \centering
    \caption{Column Partial Sums Modulo $2nk+1=681$ for the Globally Simple $H(20;17)$}
    \label{tab:H2017colsums}
    \setlength{\tabcolsep}{4pt} 
    \scriptsize
    \begin{tabular}{lc*{17}{c}} 
    \toprule
    Col 1 && 580 & 2 & 502 & 4 & 424 & 6 & \cellcolor{lightgray}649 & \cellcolor{lightgray}668 & 121 & 666 & 199 & 664 & 277 & 662 & \cellcolor{lightgray}622 & \cellcolor{lightgray}623 & \cellcolor{lightgray}0 \\ 
    Col 2 && 576 & 2 & 498 & 4 & 420 & 6 & \cellcolor{lightgray}667 & \cellcolor{lightgray}7 & 145 & 5 & 223 & 3 & 301 & 1 & \cellcolor{lightgray}97 & \cellcolor{lightgray}79 & \cellcolor{lightgray}0 \\ 
    Col 3 && 572 & 2 & 494 & 4 & 416 & 6 & \cellcolor{lightgray}665 & \cellcolor{lightgray}7 & 149 & 5 & 227 & 3 & 305 & 1 & \cellcolor{lightgray}79 & \cellcolor{lightgray}62 & \cellcolor{lightgray}0 \\ 
    Col 4 && 568 & 2 & 490 & 4 & 412 & 6 & \cellcolor{lightgray}663 & \cellcolor{lightgray}7 & 153 & 5 & 231 & 3 & 309 & 1 & \cellcolor{lightgray}96 & \cellcolor{lightgray}80 & \cellcolor{lightgray}0 \\ 
    Col 5 && 564 & 2 & 486 & 4 & 408 & 6 & \cellcolor{lightgray}661 & \cellcolor{lightgray}7 & 157 & 5 & 235 & 3 & 313 & 1 & \cellcolor{lightgray}78 & \cellcolor{lightgray}63 & \cellcolor{lightgray}0 \\ 
    Col 6 && 560 & 2 & 482 & 4 & 404 & 6 & \cellcolor{lightgray}659 & \cellcolor{lightgray}7 & 161 & 5 & 239 & 3 & 317 & 1 & \cellcolor{lightgray}95 & \cellcolor{lightgray}81 & \cellcolor{lightgray}0 \\ 
    Col 7 && 556 & 2 & 478 & 4 & 400 & 6 & \cellcolor{lightgray}657 & \cellcolor{lightgray}7 & 165 & 5 & 243 & 3 & 321 & 1 & \cellcolor{lightgray}77 & \cellcolor{lightgray}64 & \cellcolor{lightgray}0 \\ 
    Col 8 && 552 & 2 & 474 & 4 & 396 & 6 & \cellcolor{lightgray}655 & \cellcolor{lightgray}7 & 169 & 5 & 247 & 3 & 325 & 1 & \cellcolor{lightgray}94 & \cellcolor{lightgray}82 & \cellcolor{lightgray}0 \\ 
    Col 9 && 548 & 2 & 470 & 4 & 392 & 6 & \cellcolor{lightgray}653 & \cellcolor{lightgray}7 & 173 & 5 & 251 & 3 & 329 & 1 & \cellcolor{lightgray}76 & \cellcolor{lightgray}65 & \cellcolor{lightgray}0 \\ 
    Col 10 && 544 & 2 & 466 & 4 & 388 & 6 & \cellcolor{lightgray}106 & \cellcolor{lightgray}7 & 177 & 5 & 255 & 3 & 333 & 1 & \cellcolor{lightgray}93 & \cellcolor{lightgray}83 & \cellcolor{lightgray}0 \\ 
    Col 11 && 540 & 2 & 462 & 4 & 384 & 6 & \cellcolor{lightgray}650 & \cellcolor{lightgray}589 & 82 & 587 & 160 & 585 & 238 & 583 & \cellcolor{lightgray}657 & \cellcolor{lightgray}73 & \cellcolor{lightgray}0 \\ 
    Col 12 && 536 & 2 & 458 & 4 & 380 & 6 & \cellcolor{lightgray}630 & \cellcolor{lightgray}5 & 183 & 3 & 261 & 1 & 339 & 680 & \cellcolor{lightgray}83 & \cellcolor{lightgray}91 & \cellcolor{lightgray}0 \\ 
    Col 13 && 532 & 2 & 454 & 4 & 376 & 6 & \cellcolor{lightgray}632 & \cellcolor{lightgray}5 & 107 & 3 & 185 & 1 & 263 & 680 & \cellcolor{lightgray}65 & \cellcolor{lightgray}72 & \cellcolor{lightgray}0 \\ 
    Col 14 && 528 & 2 & 450 & 4 & 372 & 6 & \cellcolor{lightgray}634 & \cellcolor{lightgray}5 & 111 & 3 & 189 & 1 & 267 & 680 & \cellcolor{lightgray}84 & \cellcolor{lightgray}90 & \cellcolor{lightgray}0 \\ 
    Col 15 && 524 & 2 & 446 & 4 & 368 & 6 & \cellcolor{lightgray}636 & \cellcolor{lightgray}5 & 115 & 3 & 193 & 1 & 271 & 680 & \cellcolor{lightgray}66 & \cellcolor{lightgray}71 & \cellcolor{lightgray}0 \\ 
    Col 16 && 520 & 2 & 442 & 4 & 364 & 6 & \cellcolor{lightgray}638 & \cellcolor{lightgray}5 & 119 & 3 & 197 & 1 & 275 & 680 & \cellcolor{lightgray}85 & \cellcolor{lightgray}89 & \cellcolor{lightgray}0 \\ 
    Col 17 && 516 & 2 & 438 & 4 & 360 & 6 & \cellcolor{lightgray}640 & \cellcolor{lightgray}5 & 123 & 3 & 201 & 1 & 279 & 680 & \cellcolor{lightgray}67 & \cellcolor{lightgray}70 & \cellcolor{lightgray}0 \\ 
    Col 18 && 512 & 2 & 434 & 4 & 356 & 6 & \cellcolor{lightgray}642 & \cellcolor{lightgray}5 & 127 & 3 & 205 & 1 & 283 & 680 & \cellcolor{lightgray}86 & \cellcolor{lightgray}88 & \cellcolor{lightgray}0 \\ 
    Col 19 && 508 & 2 & 430 & 4 & 352 & 6 & \cellcolor{lightgray}644 & \cellcolor{lightgray}5 & 131 & 3 & 209 & 1 & 287 & 680 & \cellcolor{lightgray}68 & \cellcolor{lightgray}59 & \cellcolor{lightgray}0 \\ 
    Col 20 && 504 & 2 & 426 & 4 & 348 & 6 & \cellcolor{lightgray}646 & \cellcolor{lightgray}586 & 35 & 584 & 113 & 582 & 191 & 580 & \cellcolor{lightgray}678 & \cellcolor{lightgray}642 & \cellcolor{lightgray}0 \\ 
    \bottomrule
    \end{tabular}
\end{table}

\subsubsection{Column Partial Sums Are Distinct Modulo $2nk+1$}

To prove that this construction produces a globally simple array for the general case, we first state several lemmas used to derive the column partial sums for specific sets of columns.

\begin{lem}
Let $H(n;5)$ be the base array generated via the Dinitz and Wanless construction \cite{dinitz2017}. The non-empty entries of this array are exactly those detailed in \Cref{tab:Hn5colentries0}, positioned such that in each column $c$, two consecutive entries start at $(c+h-1, c)$ and three consecutive entries start at $(c-2, c)$.
\end{lem}

\begin{table}[h]
\centering
\caption{Base $H(n;5)$ Column Entries}
\label{tab:Hn5colentries0}
\setlength{\tabcolsep}{4pt} 
\footnotesize
\begin{tabular}{lcccccc}
\toprule
\textbf{Column} & \multicolumn{2}{c}{\textbf{Two Consecutive Entries}} && \multicolumn{3}{c}{\textbf{Three Consecutive Entries}} \\
\midrule
Column 1 & $-2n+2$, & $n-1$, && $-2n$, & $1$, & $3n-2$ \\
\midrule
Even $2 \le c \le h-2$ & $-n-2c+4$, & $n+2c-3$, && $5n-c/2-3$, & $-n+c$, & $-4n-c/2+2$ \\
\addlinespace
Odd $3 \le c \le h-1$ & $-n-2c+4$, & $n+2c-3$, && $4n-(c-1)/2-1$, & $-n+c$, & $-3n-(c-1)/2-1$ \\
\midrule
Column $h$ & $5n$, & $-5n+1$, && $19n/4-3$, & $-n/2$, & $-17n/4+2$ \\
\addlinespace
Column $h+1$ & $-2n+3$, & $-3n-1$, && $15n/4-1$, & $5n-3$, & $-15n/4+2$ \\
\midrule
Even $h+2 \le c \le n-2$ & $-4n+2c-1$, & $4n-2c$, && $4n+c/2-2$, & $n-c$, & $-5n+c/2+3$ \\
\addlinespace
Odd $h+3 \le c \le n-3$ & $-4n+2c-1$, & $4n-2c$, && $3n+(c-1)/2$, & $n-c$, & $-4n+(c-1)/2+2$ \\
\midrule
Column $n-1$ & $-2n-3$, & $2n+2$, && $7n/2-1$, & $-n/2+1$, & $-3n+1$ \\
\addlinespace
Column $n$ & $-2n-1$, & $-3n$, && $5n-2$, & $-2n+4$, & $2n-1$ \\
\bottomrule
\end{tabular}
\end{table}

\begin{lem}
The entries in diagonals $\mathbf{D_{2+2i}}$ and $\mathbf{D_{3+2i}}$ are given by the alternating pairs with $1 \le i \le t$:
\begin{equation*}
  -n-4ni+3-4c \quad \text{and} \quad n+4ni-1+4c.
\end{equation*}
Similarly, the entries in column $c$ for diagonals $\mathbf{D_{h+1+2i}}$ and $\mathbf{D_{h+2+2i}}$ are given by the alternating pairs:
\begin{equation*}
  n+4ni+2+4((c-h-3) \bmod{n}) \quad \text{and} \quad -n-4ni-4-4((c-h-3) \bmod{n}).
\end{equation*}
\end{lem}

\begin{lem}\label{lem:coladdedsums0_D23}
For any column $1 \le c \le n$, summing the entries of the diagonals $\{\mathbf{D}_{2+2i}, \mathbf{D}_{3+2i}\}$ top-to-bottom starting in cell $(c+3, c)$ with $1 \le i \le t$ yields the first $2t$ column partial sums $S_C(c, 1), \dots, S_C(c, 2t)$ as the sequence:
\begin{equation*}
-5n-4c+3, \ 2, \ -9n-4c+5, \ 4, \ \dots, \ -n-(4n-2)t-4c+1, \ 2t.
\end{equation*}
\end{lem}

\begin{lem}\label{lem:coladdedsums0_Dh}
For any column $1 \le c \le n$, summing the entries of the diagonals $\{\mathbf{D}_{h+1+2i}, \mathbf{D}_{h+2+2i}\}$ top-to-bottom starting in cell $(c+h+1, c)$ with $1 \le i \le t$ yields a sequence of $2t$ partial sums $S_C(c, x+1), \dots, S_C(c, x+2t)$, where $\mathcal{S} = S_C(c,x)$ is the preceding partial sum. 

For columns $1 \le c \le h+2$, the sequence is:
\begin{equation*}
\mathcal{S}+7n+4c-10, \ \mathcal{S}-2, \ \mathcal{S}+11n+4c-12, \ \mathcal{S}-4, \ \dots, \ \mathcal{S}+3n+(4n-2)t+4c-8, \ \mathcal{S}-2t.
\end{equation*}

For columns $h+3 \le c \le n$, the sequence is:
\begin{equation*}
\mathcal{S}+3n+4c-10, \ \mathcal{S}-2, \ \mathcal{S}+7n+4c-12, \ \mathcal{S}-4, \ \dots, \ \mathcal{S}-n+(4n-2)t+4c-8, \ \mathcal{S}-2t.
\end{equation*}
\end{lem}

Using these lemmas, we construct the complete partial sum sequence $S_C(c, p)$. We evaluate the first $k-3$ terms of $S_C(c, p)$ by sequentially adding: the $2t$ entries from $\mathbf{D_{2+2i}}$ and $\mathbf{D_{3+2i}}$ in \Cref{lem:coladdedsums0_D23}, the two consecutive base $H(n;5)$ entries (\cref{tab:Hn5colentries0}), and the $2t$ entries from $\mathbf{D_{h+1+2i}}$ and $\mathbf{D_{h+2+2i}}$ in \Cref{lem:coladdedsums0_Dh}. The structure of the base $H(n;5)$ array requires us to partition these sums into seven cases depending on column $c$: $c=1$, $2 \le c \le h-1$, $c=h$, $c=h+1$, $c=h+2$, $h+3 \le c \le n-1$, and $c=n$.

To see that the partial sums within each column case are distinct, observe that they partition into disjoint intervals. With the exception of columns $1$, $h+1$, and $n$, the first $k-3$ partial sums are contained in a set of partial sums bounded between $1$ and $2t+1$, a second set of partial sums bounded between $5n+2t$ and $4nt+5n-1$, and a third set of partial sums close to $2nk$. Since the individual sequences comprising these sets are strictly monotonic with respect to $i$ and occupy mutually disjoint sub-intervals, all $k-3$ partial sums are guaranteed to be distinct. Recall that $k \ge 13$ implies $t=(k-5)/4 \ge 2$, and $k=4t+5$ implies $2nk = 8nt+10n$. 

In column $1$, the $k-3$ partial sums are contained in four sets; the first set of partial sums is bounded between $2$ and $2t$, the second bounded between $6n+2t-5$ and $4nt+2n-3$, the third bounded between $2nk-4nt-n+2t-2$ and $2nk-5n$, and the fourth bounded between $2nk-2n+2t+3$ and $2nk-n+2t+2$.  

In column $h+1$, the $k-3$ partial sums are contained in three sets; the first set of partial sums is bounded between $2$ and $2t$, the second bounded between $4n+2t-4$ and $4nt-2$, and the third bounded between $2nk-4nt-n-4h+2t-2$ and $2nk-2n+2t+4$.

In column $n$, the $k-3$ partial sums are contained in three sets and a single large value; the first set of partial sums is bounded between $2$ and $2t$, the second bounded between $2n+2t-11$ and $4nt-2n-9$, the third bounded between $2nk-4nt-5n+2t+2$ and $2nk-5n+2t$, and a final large value of $2nk-2n+2t$.    

By evaluating these cumulative additions and reducing negative sums modulo $2nk+1$, we obtain the following sequences of $k-3$ partial sums for each column case:

\paragraph{Column $1$:}
After the first $2t$ partial sums, the sum is $S_C(1, 2t) \equiv 2t$. We add the two base $H(n;5)$ entries for column $1$, which are $-2n+2$ and $n-1$, yielding $S_C(1, 2t+1) = 2nk-2n+2t+3$ and $S_C(1, 2t+2) = 2nk-n+2t+2$. Finally, we list the subsequent partial sums:
\small
\begin{align*}
    & 2nk-5n, 2, \ 2nk-9n+2, \ 4, \dots, 2nk-4nt-n+2t-2, \ 2t, \ \textbf{2nk-2n+2t+3}, \\
    & \textbf{2nk-n+2t+2}, \ 6n+2t-5, \ 2nk-n+2t, \ 10n+2t-7, \ 2nk-n+2t-2, \dots, 4nt+2n-3, \ 2nk-n+2
\end{align*}
\normalsize

\paragraph{Columns $2 \le c \le h-1$:}
After the first $2t$ partial sums, the sum is $S_C(c, 2t) \equiv 2t$. We add the two base $H(n;5)$ entries for this column range, which are $-n-2c+4$ and $n+2c-3$, yielding $S_C(c, 2t+1) = 2nk-n-2c+2t+5$ and $S_C(c, 2t+2) = 2t+1$. Finally, we list the subsequent partial sums:
\small
\begin{align*}
    & 2nk-5n-4c+4, \ 2, \ 2nk-9n-4c+6, \ 4, \dots, 2nk-4nt-n-4c+2t+2, \ 2t, \ \textbf{2nk-n-2c+2t+5},\\
    & \textbf{2t+1}, \ 7n+4c+2t-9, \ 2t-1, \ 11n+4c+2t-11, \ 2t-3, \dots, 4nt+3n+4c-7, \ 1
\end{align*}
\normalsize

\paragraph{Column $h$:}
After the first $2t$ partial sums, the sum is $S_C(h, 2t) \equiv 2t$. We add the two base $H(n;5)$ entries for column $h$, which are $5n$ and $-5n+1$, yielding $S_C(h, 2t+1) = 5n+2t$ and $S_C(h, 2t+2) = 2t+1$. Finally, we list the subsequent partial sums:
\small
\begin{align*}
    & 2nk-5n-4h+4, \ 2, \ 2nk-9n-4h+6, \ 4, \dots, 2nk-4nt-n-4h+2t+2, \ 2t, \ \textbf{5n+2t},\\
    & \textbf{2t+1}, \ 7n+4h+2t-9, \ 2t-1, \ 11n+4h+2t-11, \ 2t-3, \dots, 4nt+3n+4h-7, \ 1
\end{align*}
\normalsize

\paragraph{Column $h+1$:}
After the first $2t$ partial sums, the sum is $S_C(h+1, 2t) \equiv 2t$. We add the two base $H(n;5)$ entries for column $h+1$, which are $-2n+3$ and $-3n-1$, yielding $S_C(h+1, 2t+1) = 2nk-2n+2t+4$ and $S_C(h+1, 2t+2) = 2nk-5n+2t+3$. Finally, we list the subsequent partial sums:
\small
\begin{align*}
    & 2nk-5n-4h, 2, 2nk-9n-4h+2, 4, \dots, 2nk-4nt-n-4h+2t-2, 2t, \ \textbf{2nk-2n+2t+4},\\
    & \textbf{2nk-5n+2t+3}, \ 4n+2t-4, 2nk-5n+2t+1, 8n+2t-6, 2nk-5n+2t-1, \dots, 4nt-2, 2nk-5n+3
\end{align*}
\normalsize

\paragraph{Column $h+2$:}
After the first $2t$ partial sums, the sum is $S_C(h+2, 2t) \equiv 2t$. We add the two base $H(n;5)$ entries for column $h+2$, which are $-4n+2h+3$ and $4n-2h-4$, yielding $S_C(h+2, 2t+1) = 2nk-4n+2h+2t+4$ and $S_C(h+2, 2t+2) = 2t-1$. Finally, we list the subsequent partial sums:
\small
\begin{align*}
    & 2nk-5n-4h-4, \ 2, \ 2nk-9n-4h-2, \ 4, \dots, 2nk-4nt-n-4h+2t-6, \ 2t, \ \textbf{2nk-4n+2h+2t+4},\\
    & \textbf{2t-1}, \ 7n+4h+2t-3, \ 2t-3, \ 11n+4h+2t-5, \ 2t-5, \dots, 4nt+3n+4h-1, \ 2nk
\end{align*}
\normalsize

\paragraph{Columns $h+3 \le c \le n-1$:}
After the first $2t$ partial sums, the sum is $S_C(c, 2t) \equiv 2t$. We add the two base $H(n;5)$ entries for this column range, which are $-4n+2c-1$ and $4n-2c$, yielding $S_C(c, 2t+1) = 2nk-4n+2c+2t$ and $S_C(c, 2t+2) = 2t-1$. Finally, we list the subsequent partial sums:
\small
\begin{align*}
    & 2nk-5n-4c+4, \ 2, \ 2nk-9n-4c+6, \ 4, \dots, 2nk-4nt-n-4c+2t+2, \ 2t, \ \textbf{2nk-4n+2c+2t},\\
    & \textbf{2t-1}, \ 3n+4c+2t-11, \ 2t-3, \ 7n+4c+2t-13, \ 2t-5, \dots, 4nt-n+4c-9, \ 2nk
\end{align*}
\normalsize

\paragraph{Column $n$:}
After the first $2t$ partial sums, the sum is $S_C(n, 2t) \equiv 2t$. We add the two base $H(n;5)$ entries for column $n$, which are $-2n-1$ and $-3n$, yielding $S_C(n, 2t+1) = 2nk-2n+2t$ and $S_C(n, 2t+2) = 2nk-5n+2t$. Finally, we list the subsequent partial sums:
\small
\begin{align*}
    & 2nk-9n+4, \ 2, \ 2nk-13n+6, \ 4, \dots, 2nk-4nt-5n+2t+2, \ 2t, \ \textbf{2nk-2n+2t},\\
    & \textbf{2nk-5n+2t}, \ 2n+2t-11, \ 2nk-5n+2t-2, \ 6n+2t-13, \ 2nk-5n+2t-4, \dots, 4nt-2n-9, \ 2nk-5n
\end{align*}
\normalsize

We note three properties of these calculations. First, the $k-3$ partial sums in any column are distinct and nonzero modulo $2nk+1$. Second, there are only five possible values for the partial sum $S_C(c, k-3)$ in any column: $1$, $2nk-5n$, $2nk-5n+3$, $2nk-n+2$, and $2nk$. Third, with the exception of column $h+1$, no partial sum $S_C(c, p)$ (for $1 \le p \le k-3$) in this section falls within the range $3n-1 \le S_C(c, p) \le 5n-3$, while column $h+1$ has a single partial sum in this range, $4n+2t-4$. 

Using these sequences, we calculate the final three partial sums using the values in \cref{tab:Hn5colentries0}. The results are compiled in \cref{tab:finalthreesums0}. We observe that for the general columns, with the exception of $1$, $h+1$, and $n$, the final three partial sums, $S_C(c, k-2), S_C(c, k-1)$, and $S_C(c, k)$, are either $0$ or fall strictly within the interval $[3n-1, 5n-3]$. Because this interval is avoided by the first $k-3$ partial sums, for these columns, the partial sums are distinct. For the columns $1$, $h+1$, and $n$, distinctness is similarly verified by examining the boundaries of their sets of partial sums:

Column $1$: The final nonzero partial sums, $2nk-3n+2$ and $2nk-3n+3$, fall strictly above the upper bound $2nk-5n$ of the third set and below the lower bound $2nk-2n+2t+3$ of the fourth set of partial sums.

Column $h+1$: The final nonzero partial sums, $2nk-5n/4+2$ and $15n/4-2$, fall outside the sets. Specifically, $15n/4-2$ falls strictly above the upper bound of the first set and below the lower bound of the second set, while $2nk-5n/4+2$ exceeds the third set's upper bound of $2nk-2n+2t+4$. 

Column $n$: The final nonzero partial sums, $2nk-2$ and $2nk-2n+2$, fall above the third set's upper bound $2nk-5n+2t$ and are distinct from the single large partial sum $2nk-2n+2t$ since $t \ge 2$. 

Because these final intermediate sums are bounded entirely within the empty gaps of their respective sequences, they cannot duplicate any prior term. Therefore, the columns of Heffter arrays constructed in this manner are globally simple.

\begin{table}[h]
\centering
\caption{Final Three Column Partial Sums}
\label{tab:finalthreesums0}
\setlength{\tabcolsep}{4pt} 
\footnotesize
\begin{tabular}{lccccc}
\toprule
\textbf{Column} & \textbf{$\mathbf{S_C(c, k-3)}$} && $\mathbf{S_C(c, k-2)}$ & $\mathbf{S_C(c, k-1)}$ & $\mathbf{S_C(c, k)}$ \\
\midrule
Column 1 & $2nk-n+2$ && $2nk-3n+2$ & $2nk-3n+3$ & $0$ \\
\midrule
Even $2 \le c \le h-2$ & $1$ && $5n-c/2-2$ & $4n+c/2-2$ & $0$ \\
\addlinespace
Odd $3 \le c \le h-1$ & $1$ && $4n-(c-1)/2$ & $3n+(c-1)/2+1$ & $0$ \\
\midrule
Column $h$ & $1$ && $19n/4-2$ & $17n/4-2$ & $0$ \\
\addlinespace
Column $h+1$ & $2nk-5n+3$ && $2nk-5n/4+2$ & $15n/4-2$ & $0$ \\
\midrule
Even $h+2 \le c \le n-2$ & $2nk$ && $4n+c/2-3$ & $5n-c/2-3$ & $0$ \\
\addlinespace
Odd $h+3 \le c \le n-3$ & $2nk$ && $3n+(c-1)/2-1$ & $4n-(c-1)/2-2$ & $0$ \\
\midrule
Column $n-1$ & $2nk$ && $7n/2-2$ & $3n-1$ & $0$ \\
\addlinespace
Column $n$ & $2nk-5n$ && $2nk-2$ & $2nk-2n+2$ & $0$ \\
\bottomrule
\end{tabular}
\end{table}

\subsection{Rows are simple}

To prove a row is simple, we compute its partial sums and ensure they are pairwise distinct. For each row $r \in \{1,2,\dots,n\}$, the row partial sums $S_R(r, p)$ are calculated left-to-right, starting in cell $(r, r+4)$ and wrapping around modulo $n$. As before, this choice of starting cell does not affect simplicity.

\subsubsection{Example: Globally Simple $H(20;17)$ Row Partial Sums}
We illustrate this process using the globally simple $H(20;17)$ in \Cref{tab:H2017}. \Cref{tab:H2017_Row_Sums} displays the computed partial sums for each row. For instance, the values in the first row of the table correspond to the sequence $S_R(1, 1), S_R(1, 2), \dots, S_R(1, 17)$. Shaded cells indicate partial sums derived from the base array $H(20;5)$.

\begin{table}[h]
    \centering
    \caption{Row Partial Sums Modulo $2nk+1=681$ for Globally Simple $H(20;17)$}
    \label{tab:H2017_Row_Sums}
    \setlength{\tabcolsep}{4pt}
    \scriptsize
    \begin{tabular}{lc*{17}{c}}
    \toprule
    Row 1 && 369 & 2 & 443 & 4 & 517 & 6 &\cellcolor{lightgray} 626 &\cellcolor{lightgray} 569 & 199 & 567 & 125 & 565 & 51 & 563 &\cellcolor{lightgray} 621 &\cellcolor{lightgray} 603 &\cellcolor{lightgray} 0 \\
    Row 2 && 365 & 2 & 439 & 4 & 513 & 6 &\cellcolor{lightgray} 62 &\cellcolor{lightgray} 7 & 322 & 5 & 248 & 3 & 174 & 1 &\cellcolor{lightgray} 603 &\cellcolor{lightgray} 586 &\cellcolor{lightgray} 0 \\
    Row 3 && 361 & 2 & 435 & 4 & 509 & 6 &\cellcolor{lightgray} 60 &\cellcolor{lightgray} 7 & 326 & 5 & 252 & 3 & 178 & 1 &\cellcolor{lightgray} 620 &\cellcolor{lightgray} 604 &\cellcolor{lightgray} 0 \\
    Row 4 && 357 & 2 & 431 & 4 & 505 & 607 &\cellcolor{lightgray} 659 &\cellcolor{lightgray} 608 & 250 & 606 & 176 & 604 & 102 & 1 &\cellcolor{lightgray} 602 &\cellcolor{lightgray} 587 &\cellcolor{lightgray} 0 \\
    Row 5 && 353 & 2 & 427 & 4 & 581 & 6 &\cellcolor{lightgray} 56 &\cellcolor{lightgray} 7 & 334 & 5 & 260 & 3 & 106 & 1 &\cellcolor{lightgray} 619 &\cellcolor{lightgray} 605 &\cellcolor{lightgray} 0 \\
    Row 6 && 349 & 2 & 423 & 605 & 497 & 607 &\cellcolor{lightgray} 655 &\cellcolor{lightgray} 608 & 258 & 606 & 184 & 3 & 110 & 1 &\cellcolor{lightgray} 601 &\cellcolor{lightgray} 588 &\cellcolor{lightgray} 0 \\
    Row 7 && 345 & 2 & 499 & 4 & 573 & 6 &\cellcolor{lightgray} 52 &\cellcolor{lightgray} 7 & 342 & 5 & 188 & 3 & 114 & 1 &\cellcolor{lightgray} 618 &\cellcolor{lightgray} 606 &\cellcolor{lightgray} 0 \\
    Row 8 && 341 & 603 & 415 & 605 & 489 & 607 &\cellcolor{lightgray} 651 &\cellcolor{lightgray} 608 & 266 & 5 & 192 & 3 & 118 & 1 &\cellcolor{lightgray} 600 &\cellcolor{lightgray} 589 &\cellcolor{lightgray} 0 \\
    Row 9 && 417 & 2 & 491 & 4 & 565 & 6 &\cellcolor{lightgray} 48 &\cellcolor{lightgray} 7 & 270 & 5 & 196 & 3 & 122 & 1 &\cellcolor{lightgray} 617 &\cellcolor{lightgray} 607 &\cellcolor{lightgray} 0 \\
    Row 10 && 413 & 2 & 487 & 4 & 561 & 6 &\cellcolor{lightgray} 627 &\cellcolor{lightgray} 589 & 175 & 587 & 101 & 585 & 27 & 583 &\cellcolor{lightgray} 500 &\cellcolor{lightgray} 597 &\cellcolor{lightgray} 0 \\
    Row 11 && 409 & 2 & 483 & 4 & 557 & 6 &\cellcolor{lightgray} 25 &\cellcolor{lightgray} 5 & 276 & 3 & 202 & 1 & 128 & 680 &\cellcolor{lightgray} 607 &\cellcolor{lightgray} 615 &\cellcolor{lightgray} 0 \\
    Row 12 && 405 & 2 & 479 & 4 & 553 & 6 &\cellcolor{lightgray} 27 &\cellcolor{lightgray} 5 & 280 & 3 & 206 & 1 & 132 & 680 &\cellcolor{lightgray} 589 &\cellcolor{lightgray} 596 &\cellcolor{lightgray} 0 \\
    Row 13 && 401 & 2 & 475 & 4 & 549 & 6 &\cellcolor{lightgray} 29 &\cellcolor{lightgray} 5 & 284 & 3 & 210 & 1 & 136 & 680 &\cellcolor{lightgray} 608 &\cellcolor{lightgray} 614 &\cellcolor{lightgray} 0 \\
    Row 14 && 397 & 2 & 471 & 4 & 545 & 6 &\cellcolor{lightgray} 31 &\cellcolor{lightgray} 5 & 288 & 3 & 214 & 1 & 140 & 680 &\cellcolor{lightgray} 590 &\cellcolor{lightgray} 595 &\cellcolor{lightgray} 0 \\
    Row 15 && 393 & 2 & 467 & 4 & 541 & 6 &\cellcolor{lightgray} 33 &\cellcolor{lightgray} 5 & 292 & 3 & 218 & 1 & 144 & 680 &\cellcolor{lightgray} 609 &\cellcolor{lightgray} 613 &\cellcolor{lightgray} 0 \\
    Row 16 && 389 & 2 & 463 & 4 & 537 & 6 &\cellcolor{lightgray} 35 &\cellcolor{lightgray} 5 & 296 & 3 & 222 & 1 & 148 & 680 &\cellcolor{lightgray} 591 &\cellcolor{lightgray} 594 &\cellcolor{lightgray} 0 \\
    Row 17 && 385 & 2 & 459 & 4 & 533 & 6 &\cellcolor{lightgray} 37 &\cellcolor{lightgray} 5 & 300 & 3 & 226 & 1 & 152 & 680 &\cellcolor{lightgray} 610 &\cellcolor{lightgray} 612 &\cellcolor{lightgray} 0 \\
    Row 18 && 381 & 2 & 455 & 4 & 529 & 6 &\cellcolor{lightgray} 39 &\cellcolor{lightgray} 5 & 304 & 3 & 230 & 1 & 156 & 680 &\cellcolor{lightgray} 592 &\cellcolor{lightgray} 583 &\cellcolor{lightgray} 0 \\
    Row 19 && 377 & 2 & 451 & 4 & 525 & 6 &\cellcolor{lightgray} 41 &\cellcolor{lightgray} 141 & 444 & 139 & 370 & 137 & 296 & 135 &\cellcolor{lightgray} 76 &\cellcolor{lightgray} 40 &\cellcolor{lightgray} 0 \\
    Row 20 && 373 & 2 & 447 & 4 & 521 & 6 &\cellcolor{lightgray} 588 &\cellcolor{lightgray} 551 & 177 & 549 & 103 & 547 & 29 & 545 &\cellcolor{lightgray} 584 &\cellcolor{lightgray} 585 &\cellcolor{lightgray} 0 \\
    \bottomrule
    \end{tabular}
\end{table}

\subsubsection{Row Partial Sums Are Distinct Modulo $2nk+1$}

As with the columns, we state several lemmas used to derive the row partial sums for specific sets of rows.

\begin{lem}
Let $H(n;5)$ be the base array generated via the Dinitz and Wanless construction \cite{dinitz2017}. The non-empty entries of this array are exactly those detailed in \cref{tab:Hn5rowentries0}, positioned such that in each row $r$, two consecutive entries start at $(r, r+h)$ and three consecutive entries start at $(r, r)$.
\end{lem}

\begin{table}[h]
\centering
\caption{Base $H(n;5)$ Row Entries}
\label{tab:Hn5rowentries0}
\setlength{\tabcolsep}{4pt} 
\footnotesize
\begin{tabular}{lcccccc}
\toprule
\textbf{Row} & \multicolumn{2}{c}{\textbf{Two Consecutive Entries}} && \multicolumn{3}{c}{\textbf{Three Consecutive Entries}} \\
\midrule
Row 1 & $-3n-1$, & $-3n+3$, && $3n-2$, & $-n+2$, & $4n-2$ \\
\midrule
Even $2 \le r \le h-2$ & $3n-2r$, & $-3n+2r+1$, && $-4n-r/2+2$, & $-n+r+1$, & $5n-r/2-4$ \\
\addlinespace
Odd $3 \le r \le  h-1$ & $3n-2r$, & $-3n+2r+1$, && $-3n-(r-1)/2-1$, & $-n+r+1$, & $4n-(r-1)/2-2$ \\
\midrule
Row $h$ & $-3n$, & $-2n+2$, && $-17n/4+2$, & $5n-3$, & $17n/4-1$ \\
\midrule
Odd $h+1 \le r \le n-3$ & $2r-3$, & $-2r+2$, && $-4n+(r-1)/2+2$, & $n-r-1$, & $3n+(r-1)/2+1$ \\
\addlinespace
Even $h+2 \le r \le n-4$ & $2r-3$, & $-2r+2$, && $-5n+r/2+3$, & $n-r-1$, & $4n+r/2-1$ \\
\midrule
Row $n-2$ & $2n-7$, & $-2n+6$, && $-9n/2+2$, & $-n/2+1$, & $5n-2$ \\
\addlinespace
Row $n-1$ & $2n-5$, & $5n$, && $-3n+1$, & $-2n+4$, & $-2n$ \\
\addlinespace
Row $n$ & $-5n+1$, & $-2n+3$, && $2n-1$, & $1$, & $5n-4$ \\
\bottomrule
\end{tabular}
\end{table}

\begin{lem}
The entries in row $r$ for diagonals $\mathbf{D_{h+2+2i}}$ and $\mathbf{D_{h+1+2i}}$ are given by the alternating pairs with $1 \le i \le t$:
\begin{equation*}
  -n-4ni-4-4((r-3-2i) \bmod{n}) \quad \text{and} \quad n+4ni+2+4((r-2-2i) \bmod{n}).
\end{equation*}
Similarly, the entries in diagonals $\mathbf{D_{3+2i}}$ and $\mathbf{D_{2+2i}}$ are given by the alternating pairs:
\begin{equation*}
  n+4ni+3+4((r-3-2i) \bmod{n}) \quad \text{and} \quad -n-4ni-1-4((r-2-2i) \bmod{n}).
\end{equation*}
\end{lem}

\begin{lem}\label{lem:rowaddedsums0_Dh}
For any row $1 \le r \le n$, summing the entries of the diagonals $\{\mathbf{D}_{h+2+2i}, \mathbf{D}_{h+1+2i}\}$ left-to-right starting in cell $(r, r+4)$ with $1 \le i \le t$ yields the first $2t$ row partial sums $S_R(r, 1), \dots, S_R(r, 2t)$. The sequence depends on $r$:
\begin{description}
    \item[Case 1: Rows $1 \le r \le 3$.] The sums appear as pairs $-9n-4nt+8t-4r+14+(4n-6)i$ and $2i$. 
    \item[Case 2: Even Rows $4 \le r \le (k-1)/2$.]
    As these rows contain the starting cells of diagonals $\mathbf{D_{h+2+2i}}$, $\mathbf{D_{h+1+2i}}$, $\mathbf{D_{3+2i}}$, and $\mathbf{D_{2+2i}}$, the partial sums jump once $i$ reaches the starting cell.
    \begin{itemize}
        \item For $i \le t - (r-2)/2$: The sums are $-9n-4nt+8t-4r+14+(4n-6)i$ and $2i$.
        \item For $i > t - (r-2)/2$: The sums are $-9n-4nt+8t-4r+14+(4n-6)i$ and $-4n+2i$.
    \end{itemize}
    \item[Case 3: Odd Rows $5 \le r \le (k+1)/2$.]
    These rows also contain starting cells.
    \begin{itemize}
        \item For $i < t - (r-5)/2$: The sums are $-9n-4nt+8t-4r+14+(4n-6)i$ and $2i$.
        \item For $i \ge t - (r-5)/2$: The sums are $-5n-4nt+8t-4r+14+(4n-6)i$ and $2i$.
    \end{itemize}
    \item[Case 4: Rows $(k+3)/2 \le r \le n$.]
    The sums appear as pairs $-5n-4nt+8t-4r+14+(4n-6)i$ and $2i$.
\end{description}
\end{lem}

\begin{lem}\label{lem:rowaddedsums0_D23}
For any row $1 \le r \le n$, summing the entries of the diagonals $\{\mathbf{D}_{3+2i}, \mathbf{D}_{2+2i}\}$ left-to-right starting in cell $(r, r)$ with $1 \le i \le t$ yields a sequence of $2t$ partial sums $S_R(r, x+1), \dots, S_R(r, x+2t)$, where $\mathcal{S} = S_R(r,x)$ is the preceding partial sum. The sequence depends on $r$:
\begin{description}
    \item[Case 1: Rows $1 \le r \le 3$.] The sums appear as pairs $\mathcal{S}+9n+4nt-8t+4r-15-(4n-6)i$ and $\mathcal{S}-2i$. 
    \item[Case 2: Even Rows $4 \le r \le (k-1)/2$.]
    \begin{itemize}
        \item For $i \le t - (r-2)/2$: The sums are $\mathcal{S}+9n+4nt-8t+4r-15-(4n-6)i$ and $\mathcal{S}-2i$.
        \item For $i > t - (r-2)/2$: The sums are $\mathcal{S}+9n+4nt-8t+4r-15-(4n-6)i$ and $\mathcal{S}+4n-2i$.
    \end{itemize}
    \item[Case 3: Odd Rows $5 \le r \le (k+1)/2$.]
    \begin{itemize}
        \item For $i < t - (r-5)/2$: The sums are $\mathcal{S}+9n+4nt-8t+4r-15-(4n-6)i$ and $\mathcal{S}-2i$.
        \item For $i \ge t - (r-5)/2$: The sums are $\mathcal{S}+5n+4nt-8t+4r-15-(4n-6)i$ and $\mathcal{S}-2i$.
    \end{itemize}
    \item[Case 4: Rows $(k+3)/2 \le r \le n$.]
    The sums appear as pairs $\mathcal{S}+5n+4nt-8t+4r-15-(4n-6)i$ and $\mathcal{S}-2i$.
\end{description}
\end{lem}

Using these lemmas, we construct the complete partial sum sequence $S_R(r, p)$. We evaluate the first $k-3$ terms of $S_R(r, p)$ by sequentially adding: the $2t$ entries from diagonals $\mathbf{D_{h+2+2i}}$ and $\mathbf{D_{h+1+2i}}$ in \Cref{lem:rowaddedsums0_Dh}, the two consecutive base $H(n;5)$ entries (\cref{tab:Hn5rowentries0}), and the $2t$ entries from diagonals $\mathbf{D_{3+2i}}$ and $\mathbf{D_{2+2i}}$ in \Cref{lem:rowaddedsums0_D23}. The structure of the base $H(n;5)$ array requires us to partition these sums into nine cases depending on row $r$: $r=1$, $r \in \{2,3\}$, even and odd rows from $4 \le r \le (k+1)/2$, $(k+3)/2 \le r \le h-1$ (which is empty when $k=n-3$), $r=h$, $h+1 \le r \le n-2$, $r=n-1$, and $r=n$.

To see that the partial sums within each row case are distinct, observe that, with the exception of row $n-1$, they partition into disjoint intervals: a set of partial sums bounded above by $2t+1$, a second set of partial sums bounded between $n+2t-3$ and $4nt+5n+2t-4$, and a third set of partial sums close to $2nk$. In row $n-1$ we have a set of partial sums bounded above by $2t$, a second set of partial sums bounded between $2n+2t-5$ and $4nt+12n-6t-18$, and a third set of partial sums close to $2nk$. Since the individual sequences comprising these sets are strictly monotonic with respect to $i$ and occupy mutually disjoint sub-intervals, all $k-3$ partial sums are guaranteed to be distinct. Recall that $k \ge 13$ implies $t=(k-5)/4 \ge 2$, and $k=4t+5$ implies $2nk = 8nt+10n$.

By evaluating these cumulative additions and reducing negative sums modulo $2nk+1$, we obtain the following sequences of $k-3$ partial sums for each row case:

\paragraph{Row $1$:} After the first $2t$ partial sums, the sum is $S_R(1, 2t) \equiv 2t$. We add the two base $H(n;5)$ entries for row $1$, which are $-3n-1$ and $-3n+3$, yielding $S_R(1, 2t+1) = 2nk-3n+2t$ and $S_R(1, 2t+2) = 2nk-6n+2t+3$. Finally, we list the subsequent partial sums:
        \small
    \begin{align*}
    & 2nk-5n-4nt+8t+5, \ 2, \ 2nk-n-4nt+8t-1, \ 4, \dots, 2nk-9n+2t+11, \ 2t, \ \textbf{2nk-3n+2t}, \\
    & \textbf{2nk-6n+2t+3}, \ 4nt-n-6t-3, \ 2nk-6n+2t+1, \dots, 3n-9, \ 2nk-6n+3
    \end{align*}
    \normalsize

\paragraph{Rows $2$ and $3$:} After the first $2t$ partial sums, the sum is $S_R(r, 2t) \equiv 2t$. We add the two base $H(n;5)$ entries for this row range, which are $3n-2r$ and $-3n+2r+1$, yielding $S_R(r, 2t+1) = 3n+2t-2r$ and $S_R(r, 2t+2) = 2t+1$. Finally, we list the subsequent partial sums:
\small
    \begin{align*}
    & 2nk-5n-4nt+8t-4r+9, \ 2, \ 2nk-n-4nt+8t-4r+3, \ 4, \dots, 2nk-9n+2t-4r+15, \ 2t, \\
    & \textbf{3n+2t-2r}, \ \textbf{2t+1}, \ 4nt+5n-6t+4r-8, \ 2t-1, \dots, 9n+4r-14, \ 1
    \end{align*}
    \normalsize

Similar to the previous construction, rows $4$ through $(k+1)/2$ contain the starting cells of diagonals $\mathbf{D_{h+2+2i}}$, $\mathbf{D_{h+1+2i}}$, $\mathbf{D_{3+2i}}$, and $\mathbf{D_{2+2i}}$. This causes the sequence of partial sums in these rows to jump when $i$ reaches the starting cell index. As detailed in \Cref{lem:rowaddedsums0_Dh} and \Cref{lem:rowaddedsums0_D23}, the partial sums will take one of two values depending on whether the starting cell has been reached. When listing the partial sums in these rows, we again use set notation to represent these two possible values. Because each of the potential sums in these sets is strictly distinct from the rest of the partial sums in the row, this notation allows us to establish the uniqueness of the sequence without the need for piecewise cases.
    
\paragraph{Even Rows $4 \le r \le (k-1)/2$:} After the first $2t$ partial sums, the sum is $S_R(r, 2t) = -4n+2t \equiv 2nk-4n+2t+1 \pmod{2nk+1}$. We add the two base $H(n;5)$ entries, $3n-2r$ and $-3n+2r+1$, yielding $S_R(r, 2t+1) = 2nk-n+2t-2r+1$ and $S_R(r, 2t+2) = 2nk-4n+2t+2$. For partial sums that vary depending on $i$ and the row, as noted in \Cref{lem:rowaddedsums0_Dh} and \Cref{lem:rowaddedsums0_D23}, we use set notation to list the possible values. Because every potential value in these sets is strictly distinct from the rest of the entries, this sufficiently proves the uniqueness of the sequence without the need for piecewise cases:
\small
    \begin{align*}
    & 2nk-5n-4nt+8t-4r+9, \ \{2, 2nk-4n+3\}, \dots, 2nk-9n+2t-4r+15, \ 2nk-4n+2t+1, \\
    & \textbf{2nk-n+2t-2r+1}, \ \textbf{2nk-4n+2t+2}, \ 4nt+n-6t+4r-8, \ \{2nk-4n+2t, 2t-1\}, \dots, 5n+4r-14, \ 1 
    \end{align*}
    \normalsize
     
\paragraph{Odd Rows $5 \le r \le (k+1)/2$:} After the first $2t$ partial sums, the sum is $S_R(r, 2t) \equiv 2t$. We add the two base $H(n;5)$ entries for this row range, which are $3n-2r$ and $-3n+2r+1$, yielding $S_R(r, 2t+1) = 3n+2t-2r$ and $S_R(r, 2t+2) = 2t+1$. As with the even rows, we use set notation to represent the partial sums that vary depending on $i$ and the row, which sufficiently demonstrates the distinctness of the sequence of partial sums:
\small
    \begin{align*}
    & \{2nk-5n-4nt+8t-4r+9, 2nk-n-4nt+8t-4r+9\}, \ 2, \dots, 2nk-5n+2t-4r+15, \ 2t, \\
    & \textbf{3n+2t-2r}, \ \textbf{2t+1}, \ \{4nt+n-6t+4r-8, 4nt+5n-6t+4r-8\}, \ 2t-1, \dots, 5n+4r-14, \ 1
    \end{align*}
    \normalsize
    
\paragraph{Rows $(k+3)/2 \le r \le h-1$:} After the first $2t$ partial sums, the sum is $S_R(r, 2t) \equiv 2t$. We add the two base $H(n;5)$ entries for this row range, which are $3n-2r$ and $-3n+2r+1$ (This section is empty when $k=n-3$.), yielding $S_R(r, 2t+1) = 3n+2t-2r$ and $S_R(r, 2t+2) = 2t+1$. Finally, we list the subsequent partial sums:
\small
    \begin{align*}
    & 2nk-n-4nt+8t-4r+9, \ 2, \ 2nk+3n-4nt+8t-4r+3, \ 4, \dots, 2nk-5n+2t-4r+15, \ 2t, \\
    & \textbf{3n+2t-2r}, \ \textbf{2t+1}, \ 4nt+n-6t+4r-8, \ 2t-1, \dots, 5n+4r-14, \ 1
    \end{align*}
    \normalsize
    
\paragraph{Row $h$:} The sum is $S_R(h, 2t) \equiv 2t$. We add the two base entries for row $h$, which are $-3n$ and $-2n+2$, yielding $S_R(h, 2t+1) = 2nk-3n+2t+1$ and $S_R(h, 2t+2) = 2nk-5n+2t+3$. Finally, we list the subsequent partial sums:
\small
    \begin{align*}
    & 2nk-3n-4nt+8t+9, \ 2, \ 2nk+n-4nt+8t+3, \ 4, \dots, 2nk-7n+2t+15, \ 2t, \\
    & \textbf{2nk-3n+2t+1}, \ \textbf{2nk-5n+2t+3}, \ 4nt-2n-6t-7, \ 2nk-5n+2t+1, \dots, 2n-13, \ 2nk-5n+3
    \end{align*}
    \normalsize
    
\paragraph{Rows $h+1 \le r \le n-2$:} From the sum $S_R(r, 2t) \equiv 2t$, we add the base entries $2r-3$ and $-2r+2$, yielding $S_R(r, 2t+1) = 2t+2r-3$ and $S_R(r, 2t+2) = 2t-1$. Finally, we list the subsequent partial sums:
\small
    \begin{align*}
    & 2nk-n-4nt+8t-4r+9, \ 2, \ 2nk+3n-4nt+8t-4r+3, \ 4, \dots, 2nk-5n+2t-4r+15, \ 2t, \\
    & \textbf{2t+2r-3}, \ \textbf{2t-1}, \ 4nt+n-6t+4r-10, \ 2t-3, \dots, 5n+4r-16, \ 2nk
    \end{align*}
    \normalsize
    
\paragraph{Row $n-1$:} From the sum $S_R(n-1, 2t) \equiv 2t$, we add the base entries $2n-5$ and $5n$, yielding $S_R(n-1, 2t+1) = 2n+2t-5$ and $S_R(n-1, 2t+2) = 7n+2t-5$. Finally, we list the subsequent partial sums:
\small
    \begin{align*}
    & 2nk-5n-4nt+8t+13, \ 2, \ 2nk-n-4nt+8t+7, \ 4, \dots, 2nk-9n+2t+19, \ 2t, \\
    & \textbf{2n+2t-5}, \ \textbf{7n+2t-5}, \ 4nt+12n-6t-18, \ 7n+2t-7, \dots, 16n-24, \ 7n-5
    \end{align*}
    \normalsize
    
\paragraph{Row $n$:} From the sum $S_R(n, 2t) \equiv 2t$, we add the base entries for row $n$, which are $-5n+1$ and $-2n+3$, yielding $S_R(n, 2t+1) = 2nk-5n+2t+2$ and $S_R(n, 2t+2) = 2nk-7n+2t+5$. Finally, we list the subsequent partial sums:
\small
    \begin{align*}
    & 2nk-5n-4nt+8t+9, \ 2, \ 2nk-n-4nt+8t+3, \ 4, \dots, 2nk-9n+2t+15, \ 2t, \\
    & \textbf{2nk-5n+2t+2}, \ \textbf{2nk-7n+2t+5}, \ 4nt-2n-6t-5, \ 2nk-7n+2t+3, \dots, 2n-11, \ 2nk-7n+5
    \end{align*}
    \normalsize

We note two properties of these calculations that allow us to easily verify the distinctness of the partial sums. First, the first $k-3$ partial sums in any row are distinct and nonzero modulo $2nk+1$. Second, there are only six possible values for the partial sum $S_R(r, k-3)$ in any row: $1$, $7n-5$, $2nk-7n+5$, $2nk-6n+3$, $2nk-5n+3$, and $2nk$. 

Using these properties, we calculate the final three partial sums using the values in \cref{tab:Hn5rowentries0}, yielding the results compiled in \cref{tab:finalthreerowsums0}. In order to show the final three partial sums, $S_R(r, k-2), S_R(r, k-1)$, and $S_R(r, k)$, are distinct from the $k-3$ sums calculated previously, we will need to consider sets of rows. 

Row $1$: The final nonzero partial sums $2nk-3n+1$ and $2nk-4n+3$ fall between the second largest, $2nk-6n+2t+3$ and largest, $2nk-3n+2t$ of the $k-3$ partial sums.

Row $2$: The final nonzero partial sums $2nk-4n+3$ and $2nk-5n+6$ fall above the largest $2nk-9n+2t+7$ of the $k-3$ partial sums.

Odd Rows $3 \le r \le h-1$: The final nonzero partial sums $2nk-3n-(r-1)/2+1$ and $2nk-4n+(r-1)/2+3$ fall above the largest $2nk-5n+2t-4r+15$ of the $k-3$ partial sums.

Even Rows $4 \le r \le h-2$: The final nonzero partial sums $2nk-4n-r/2+4$ and $2nk-5n+r/2+5$ fall between $2nk-9n+2t-4r+15$ and the next largest partial sum for any of these even rows which is $2nk-4n+3$. 

Row $h$: The nonzero partial sum $2nk-37n/4+5$ falls between the consecutive partial sums $2nk+n-4nt+8t+3$ and $2nk+5n-4nt+8t-3$, while the final nonzero partial sum $2nk-17n/4+2$ falls between the consecutive partial sums $2nk-3n+2t+1$ and $2nk-5n+2t+3$.

Rows $h+1 \le r \le n-2$: The final nonzero partial sums in these rows fall between the second largest and largest of the first $k-3$ partial sums. 

Row $n-1$: The final nonzero partial sum $2n$ falls between the consecutive partial sums $2t$ and $2n+2t-5$, while the nonzero partial sum $4n-4$ falls between the consecutive partial sums $2n+2t-5$ and $7n-5$.  

Row $n$: The final nonzero partial sums $2nk-5n+4$ and $2nk-5n+5$ between the consecutive partial sums $2nk-7n+2t+5$ and $2nk-5n+2t+2$ of the first $k-3$ partial sums.  

Therefore, all partial sums within every row are distinct. This, taken together with the fact that the columns are simple proves that the Heffter arrays constructed in this manner are globally simple.

\begin{table}[h]
\centering
\caption{Final Three Row Partial Sums}
\label{tab:finalthreerowsums0}
\setlength{\tabcolsep}{4pt} 
\footnotesize
\begin{tabular}{lccccc}
\toprule
\textbf{Row} & \textbf{$\mathbf{S_R(r, k-3)}$} && $\mathbf{S_R(r, k-2)}$ & $\mathbf{S_R(r, k-1)}$ & $\mathbf{S_R(r, k)}$ \\
\midrule
Row 1 & $2nk-6n+3$ && $2nk-3n+1$ & $2nk-4n+3$ & $0$ \\
\midrule
Even $2 \le r \le h-2$ & $1$ && $2nk-4n-r/2+4$ & $2nk-5n+r/2+5$ & $0$ \\
\addlinespace
Odd $3 \le r \le h-1$ & $1$ && $2nk-3n-(r-1)/2+1$ & $2nk-4n+(r-1)/2+3$ & $0$ \\
\midrule
Row $h$ & $2nk-5n+3$ && $2nk-37n/4+5$ & $2nk-17n/4+2$ & $0$ \\
\midrule
Odd $h+1 \le r \le n-3$ & $2nk$ && $2nk-4n+(r-1)/2+2$ & $2nk-3n-(r-1)/2$ & $0$ \\
\addlinespace
Even $h+2 \le r \le n-4$ & $2nk$ && $2nk-5n+r/2+3$ & $2nk-4n-r/2+2$ & $0$ \\
\midrule
Row $n-2$ & $2nk$ && $2nk-9n/2+2$ & $2nk-5n+3$ & $0$ \\
\addlinespace
Row $n-1$ & $7n-5$ && $4n-4$ & $2n$ & $0$ \\
\addlinespace
Row $n$ & $2nk-7n+5$ && $2nk-5n+4$ & $2nk-5n+5$ & $0$ \\
\bottomrule
\end{tabular}
\end{table}

\section{Concluding Remarks and Open Problems}

In this paper, we addressed a significant gap in the literature by explicitly constructing globally simple square Heffter arrays $H(n;k)$ for the previously open cases where $k \equiv 1 \pmod{4}$. By leveraging a base $H(n;5)$ array and systematically augmenting it with sets of four diagonals, we proved that the resulting row and column partial sums are distinct modulo $2nk+1$. Building upon these results, we identify four primary avenues for future research.

The first avenue directly extends our existence results to the two remaining open cases for globally simple square Heffter arrays. Specifically, our next objective is to adapt this framework to establish the existence of globally simple Heffter arrays when $n \equiv 0,2 \pmod{4}$ and $k \equiv 2 \pmod{4}$.

A second direction for future work concerns the biembeddings of the arrays constructed here. While this paper resolves the globally simple property for the $k \equiv 1 \pmod{4}$ arrays, the question of their biembeddings remains partially open. For the specific sub-cases where $n \equiv 3 \pmod{4}$ with $\gcd(n,k-2) > 1$, compatible orderings may still exist; successfully constructing them remains an engaging open problem.

Third, researchers could shift focus from square to rectangular Heffter arrays. It is interesting to note that there are currently relatively few constructions that are both simple and compatible in the rectangular case, and the known results primarily concern arrays with no empty cells. For instance, while there are simple and compatible constructions for all $H(3,n)$, results for simple $H(5,n)$ are currently restricted to $n \le 100$ \cite{dinitz2_2017}. Expanding this catalog represents a substantial open challenge.

Finally, a fourth broad avenue encompasses variants and generalizations of Heffter arrays. These include weak Heffter arrays, $\lambda$-fold Heffter arrays, relative Heffter arrays, non-zero sum Heffter arrays, Archdeacon arrays, and Heffter spaces \cite{archdeacon2015, buratti2024, buratti2025, costa2_2020, costa3_2020, costa2021, costa2022, costa2024, johnson2025}. Exploring connections to decompositions and embeddings within these generalized structures provides a highly fruitful area of inquiry.

\paragraph{Acknowledgement:} The authors would like to acknowledge support from TÜBİTAK 1001 124F360.

\bibliography{references}

\end{document}